\documentclass[oneside]{amsart}
\usepackage{amssymb}
\usepackage{ifthen}

\newcommand{\formatswitch}{preprint}

\newcommand{\tref}[1]{(\ref{#1})}

\ifthenelse{\equal{\formatswitch}{big}}{
\setlength{\textwidth}{432pt}
\setlength{\oddsidemargin}{18.8775pt}

\setcounter{tocdepth}{1}
}{
\ifthenelse{\equal{\formatswitch}{paper}}{

\setcounter{tocdepth}{1}
}{
\setcounter{tocdepth}{2}
}
}

\newcommand{\claimenum}{\renewcommand{\theenumi}{\alph{enumi}}
 \renewcommand{\labelenumi}{\textit{(\theenumi)}}
 \renewcommand{\theenumii}{\roman{enumii}}
 \renewcommand{\labelenumii}{\textit{(\theenumii)}}
 \begin{enumerate}}
\newcommand{\claimenumend}{\end{enumerate}}

\newtheorem{dummy}{realdumb}[section]

\newtheorem{lemma}[dummy]{Lemma}
\newtheorem{prop}[dummy]{Proposition}
{\theoremstyle{definition} }

{\theoremstyle{definition} }
\newtheorem{cor}{Corollary}[dummy]

\newcommand{\strutdepth}{\dp\strutbox}
\newcommand{\marginalnote}[1]
   {\strut\vadjust{\kern-\strutdepth\domarginalnote{#1}}}
\newcommand{\domarginalnote}[1]{\vtop to \strutdepth{
  \baselineskip\strutdepth
   \vss\llap{ #1\ \ }\null}}  %get sevenpoint font here
\newcounter{showlabelflag}
\newcounter{makelabelflag}
\newcommand{\showlabels}{\setcounter{showlabelflag}{1}}

\newcommand{\makelabels}{\setcounter{makelabelflag}{1}}
\newcommand{\hidelabels}{\setcounter{showlabelflag}{2}}
\newcommand{\mylabel}[1]{
  \ifthenelse{\value{makelabelflag}=1}
    {\label{#1}}{}
  \ifthenelse{\value{showlabelflag}=1}
    {\marginpar{#1}}{}\relax}

\newcommand{\Z}{{\mathbf Z}}
\newcommand{\N}{{\mathbf N}}
\newcommand{\sub}{\subseteq}

\ifthenelse{\equal{\formatswitch}{draft}}{
\showlabels
}{\relax}

\newcommand{\ZS}{Zappa-Sz\'ep }

\newcommand{\zapprod}{\mathbin{\bowtie}}

\newcommand{\mymargin}[1]{
  \ifthenelse{\value{showlabelflag}=1}
    {\marginpar{#1}}{}\relax}

\newcommand{\LeftAct}[2]{#1 \cdot #2}

\newcounter{keepitemnum}
\newcommand{\keepitem}{\setcounter{keepitemnum}{\value{enumi}}}
\newcommand{\useitem}{\setcounter{enumi}{\value{keepitemnum}}}
\newcounter{keepitemnumm}

\begin{document}

\bibliographystyle{amsplain}
\begin{center}{\bfseries On the \ZS Product\footnote{AMS
Classification (2000): primary 20N02, secondary 20L05, 20M10}}\end{center}
\vspace{3pt}
\begin{center}{MATTHEW G. BRIN}\end{center}
\vspace{4pt}
\vspace{3pt}
%\begin{center}\today\end{center}
\begin{center}March 3, 2003\end{center}

%\CompileMatrices

%\makelabels has to come AFTER \tableofcontents 
%to prevent the labels in the section headers from 
%being processed twice

\makelabels
%\showlabels
\hidelabels

%**end of header

\section{Introduction}

The semidirect product of two groups is a natural generalization of
the direct product of two groups in that the requirement that both
factors be normal in the product is replaced by the weaker
requirement that only one of the factors be normal in the product.
The \ZS product of two groups is a natural generalization of the
semidirect product of two groups in that neither factor is required
to be normal.

The \ZS product was developed by G. Zappa in \cite{zappa:gpprod}.
Variations and generalizations were studied in the setting of groups
by R\'edei in \cite{redei:genprod} and Casadio in
\cite{casadio:prods}.  The products were used to discover properties
of groups by R\'edei, Sz\'ep and Tibiletti in numerous papers.
Sz\'ep independently discovered the relations of the product and
used them to study structural properties of groups (e.g., normal
subgroups \cite{szep:normalsubgps}) and also initiated the study of
similar products in settings other than groups in
\cite{szep:ringprod} and \cite{szep:genprod}.  The terminology
{\itshape \ZS product} was suggested by Zappa.  See also Page 674 of
\cite{Huppert:endliche-gruppen-I}.

For the last several years, the author has been investigating and
extending a family of closely related groups known as Thompson's
groups.  See \cite{CFP} for background on these groups.  It was
found that the structure of these groups as groups of fractions of
monoids that in most cases were \ZS products of simpler monoids was
the key to successful analysis.  

The author's investigations were motivated in the first place by an
intimate relationship that Thompson's groups have with certain
categories.  The intimate relationship is helped by the fact that
the \ZS product works fantastically well with categories.  The
(partial) multiplication on a cateogry is the composition of its
morphisms.  

It turns out that the \ZS product has the remarkable property that
it seems to require no hypotheses at all.  Group-like properties
such as associativity, fullness of the multiplication, identities,
and inverses can be assumed or removed at random and the \ZS product
can still be discussed and used at some level.

This paper records the author's observations on the \ZS product that
are needed for his investigations.  Properties needed to work with
categories, monoids, and groups of fractions of monoids include
associativity, existence of certain identities and inverses,
cancellativity, and the existence of common right multiples.  The
behavior of the \ZS product in connection with these properties is
studied in this paper.  Definitions for these properties and the
property in the next paragraph are given in the next section.

Thompson's groups are unusual in that they have strong finiteness
properties (most are finitely presented) and some interesting
geometric properties.  The finiteness and geometric properties
include the ability to act nicely on certain ``non-positively
curved'' metric spaces and spaces that are finite in each dimension
modulo the action.  Crucial to the construction of these spaces is
the existence of least common left multiples in the generating
monoids.  Thus the behavior of the \ZS product in regards to
presentations and least common left multiples is also studied in
this paper.

Some of this paper is expository.  The \ZS product is ``done over''
from scratch in excruciating generality.  It is hoped that the
increased generality and the previously unexplored behavior of the
\ZS product with respect to certain aspects of algebraic structures
will be of interest to others.

The next section reviews all of the terms used in the paper.  It is
a long section since the behavior of the \ZS product is studied with
respect to several properties and constructs.  The section after
that defines the \ZS product and studies its properties.

The author would like to thank B. Brewster, D. Karagueuzian,
J. Sz\'ep, and G. Zappa for helpful comments and communications.

\section{Basics and definitions}

We use a large number of structures in this paper.  This purpose of
this section is to establish terminology.  We work with
presentations of semigroups, monoids and groups and we write out
carefully what is meant in each case.  A technique involving
confluent and terminating relations is a powerful tool that we use,
and we write out most of its details.

\subsection{Sets with multiplication} A triple \((A, D, m)\) is a
{\itshape set with partial multiplication} if \(A\) is a set,
\(D\sub A\times A\) and \(m:D\rightarrow A\) is a function.  There
are no other restrictions.  We write \(ab\) for \(m(a,b)\).  We call
\(D\) the {\itshape domain of the multiplication} and also regard it
as a relation on \(A\).  We write \(aDb\) if \((a,b)\) is in \(D\)
and we will often say that \(ab\) {\itshape exists} if \(aDb\).  The
morphisms of a (small) category form a set with a partial
multiplication where the multiplication is composition.

A {\itshape homomorphism of sets with partial multiplication}
\(f:(A,D,m)\rightarrow (B,E,n)\) is a function \(f:A\rightarrow B\)
so that \(f(D)\sub E\) and for all \((a,b)\in D\) we have
\(f(ab)=f(a)f(b)\).  

If \((A,D,m)\) is a set with a partial multiplication and \(B\sub
A\), then we say that \(B\) is {\itshape closed} under the
multiplication if \(b_1b_2\in B\) for every \((b_1,b_2)\in D\cap
(B\times B)\).

\subsection{Multiplication properties} Let \((A, D, m)\) be a set
with a partial multiplication.  We say that the partial
multiplication is {\itshape right associative} if \(aDb\) and
\((ab)Dc\) together imply that \(bDc\) and \(aD(bc)\) hold and
\((ab)c=a(bc)\).  The multiplication is {\itshape left associative}
if \(bDc\) and \(aD(bc)\) together imply that \(aDb\) and \((ab)Dc\)
hold and \((ab)c=a(bc)\).  We say that the partial multiplication is
{\itshape (fully) associative} if it is both left and right
associative.  The multiplication is {\itshape categorical} if it is
fully associative and \(aDb\) and \(bDc\) together imply that
\(aD(bc)\) and \((ab)Dc\).  The multiplication is {\itshape full} if
\(D = A\times A\).  The morphisms of a small category under
composition form a categorical multiplication which in general is
not full.

A set with a full, associative multiplication is a {\itshape
semigroup}.

Let \((A, D, m)\) be a set with a partial multiplication and let
\(x\) be in \(A\).  A {\itshape right identity} for \(x\) is an
\(a\) in \(A\) for which \(xDa\) and \(xa=x\).  We define {\itshape
left identity} for \(x\) symmetrically.  We say that \(a\in A\) is a
{\itshape right identity} for \(A\) if it is a right identity for
every \(x\) with \(xDa\) and there is at least one \(x\) with
\(xDa\).  We define {\itshape left identity} for \(A\)
symmetrically.  We say that \(a\in A\) is a {\itshape full identity}
if it is both a right and left identity for \(A\).  An element \(e\)
of a partial multiplication \((A, D, m)\) is a {\itshape global
identity} if, for all \(a\in A\), we have \(eDa\) and \(aDe\) and
\(ae=ea=a\).  Identities of various types are far from unique.
However, we have the following.

\begin{lemma}\mylabel{UniqueId} In a right (left) associative
partial multiplication, there can be only one right (left) identity
for a given \(x\) that is also a full identity.  \end{lemma}

\begin{proof} We give the proof using the word ``right.''  If \(a\)
and \(b\) are right identities for \(x\), then \((xa)b=xb=x\) and
right associativity requires that \(ab\) must exist.  If both \(a\)
and \(b\) are full identities, then \(ab\) equals both \(a\) and
\(b\).  \end{proof}

We say that a set \((A, D, m)\) with partial multiplication
{\itshape has right (left) identities} if every element has a right
(left) identity that is also a right (left) identity for \(A\).
This is stronger than requiring that each element has a right (left)
identity.  We say that a partial multiplication has {\itshape full
identities} if it has both right and left identities.  The next
lemma shows that having full identies really means that there are
lots of full identities.

\begin{lemma}\mylabel{FullId} If a set \((A, D, m)\) with partial
multipliation has full identities, then each right (left) identity
\(a\) for \(A\) is a full identity for \(A\) and satisfies \(aDa\)
(and is thus an idempotent).  \end{lemma}

\begin{proof} It suffices to consider a right identity \(a\) for
\(A\).  There is a left identity \(b\) for \(a\) that is a left
identity for \(A\).  Now \(ba\) equals both \(a\) and \(b\), so
\(a\) is both a left and a right identity for \(A\).  Now that
\(a=b\), we have \(aDa\) and \(a^2=a\).  \end{proof}

Since there are many types of identities, it is clear that we can
invent even more types of inverses.  We will define exactly one type
since it is the only one that we will prove anything about.  A set
\(M\) with partial multiplication has {\itshape left inverses with
respect to right identities} if there is an \(x\) making \(xa=b\)
whenever \(b\) is a right identity for \(M\) and \(ab\) is defined.

We say that a partial multiplication is {\itshape left cancellative}
if \(ab=ac\) always implies \(b=c\), and say it is {\itshape right
cancellative} if \(ac=bc\) always implies that \(a=b\).  The
multiplication is said to be {\itshape cancellative} if it is both
right and left cancellative.  A partial multiplication is {\itshape
strongly left cancellative} if it is left cancellative and if
\(ab=a\) implies that \(b\) is a global identity.  The reader can
define strongly right cancellative.

\begin{lemma}\mylabel{GetStrLftCanc} Let \(S\) be a cancellative
semigroup.  Then (1) \(S\) is strongly left cancellative and
strongly right cancellative and (2) if \(S\) has a global identity
1, and \(ab=1\), then \(ba=1\).  \end{lemma}

\begin{proof} If \(ab=a\), then the chain of implications \[ab = a
\quad \Rightarrow \quad abc = ac \quad \Rightarrow \quad bc = c
\quad \Rightarrow \quad cbc = cc \quad \Rightarrow \quad cb = c \]
where \(c\) is arbitrary shows that \(b\) is a global identity.
Strong right cancellativity is argued the same way.  If \(ab=1\),
then \(aba=a\) and \(ba\) must also be a global identity by strong
left cancellativity.  \end{proof}

We say that \(a\) is a {\itshape unit} in a set \(M\) with partial
multiplication if there is a \(b\) so that both \(ab\) and \(ba\)
are a global identity for \(M\).  Of course, this makes \(b\) a unit
as well.

A {\itshape monoid} is a semigroup with a global identity.

Let \(a\) and \(b\) be elements of a semigroup for which there are
\(p\) and \(q\) so that \(ap=bq\).  We call \(ap=bq\) a {\itshape
common right multiple} of \(a\) and \(b\).  A semigroup is said to
have {\itshape common right multiples} if every pair of elements has
a common right multiple.  A symmetric definition gives common left
multiples.  We asymmetrically define the following.  The notion is
needed in certain constructions based on the Thompson groups.  A
semigroup has {\itshape least common left multiples} if for every
pair \(a\) and \(b\) with a common left multiple there are \(p\) and
\(q\) so that \(pa=qb\) and if \(m\) is another common left multiple
of \(a\) and \(b\), then there is a \(k\) so that \(m=kpa=kbq\).
Note that the definition of least common left multiples is worded so
that a semigroup with least common left multiples need not have
common left multiples.

If \(x\) and \(y\) are in a set with multiplication, then a
{\itshape common left factor} of \(x\) and \(y\) is an \(f\) so that
\(x=fa\) and \(y=fb\) for some \(a\) and \(b\).

The next lemma is straightforward.

\begin{lemma}\mylabel{MaxesUnique} Let \(S\) be a cancellative
semigroup.  Then the following hold.  {\claimenum \item Let \(l\) be
a least common left multiple of \((a,b)\).  Then \(l'\) is a least
common left multiple of \((a,b)\) if and only if \(l'=ul\) for some
unit \(u\) of \(S\).  \item Let \((a,b)\) have a least common left
multiple.  Then \(xa=yb\) is a least common left multiple of
\((a,b)\) if and only if all common left factors of \((x,y)\) are
units of \(S\).  \claimenumend} \end{lemma}

\subsection{Categories}\mylabel{CatBaseSect} We use the
standard notions of category theory such as can be found in
\cite{MacLane:categories}.

A {\itshape groupoid} is a category in which every morphisms is an
isomorphism.  A groupoid with one object forms a {\itshape group}
(with the emphasis being on the set of morphisms and its operation
of composition).  A category with one object forms a {\itshape
monoid} (with the emphasis being on the set of morphisms and its
operation of composition).

A groupoid is {\itshape connected} if every pair of objects has a
morphism between them.

If \(G\) is a groupoid and \(x\) is an object of \(G\), then \(G_x\)
will denote the set of morphisms that have \(x\) as both domain and
codomain.  Composition makes \(G_x\) into a group.  If \(G\) is
connected, then all the \(G_x\) are isomorphic.  However, there is
no canonical isomorphism between any two of them.

If \(G\) is an arbitrary category, then \(G_x\) will still denote
the set of morphisms that have \(x\) as both domain and codomain,
but now \(G_x\) is just a monoid.  Even if \(G\) is connected
(again, any two objects have a morphism between them) it is not
necessarily the case that all the \(G_x\) are isomorphic.

\subsection{Relations} Relations are subsets of cross products and
as such can be operated on by the usual set operations such as union
and intersection.  Relations with properties that are preserved
under arbitrary intersection (such as equivalence relations) can be
generated by (sets of) other relations.  Thus we can talk of the
equivalence relation generated by one or more relations on a set
\(A\) as the smallest equivalence relation on \(A\) that contains
all the given relations.  It is the intersection of all the
equivalence relations containing the given relations and it is
provably the set of all pairs \((a,b)\) for which there is a chain
of elements \(a=a_0,a_1,\ldots,a_k=b\) so that, for
\(i\in\{0,1,\ldots,k-1\}\), one of \((a_i,a_{i+1})\) or
\((a_{i+1},a_i)\) is in one of the given relations or
\(a_i=a_{i+1}\).  We can also talk about the symmetric closure
(smallest containing symmetric relation) or transitive closure,
etc., of a set of relations.  The equivalence relation generated by
a set of relations is the reflexive, symmetric, transitive closure.

The empty relation is allowed, as is the empty function.

If \(\rightarrow\) is a binary relation on a set \(A\), then we
write \(a\rightarrow b\) for \((a,b)\in \rightarrow\).  The inverse
relation \(\leftarrow\) is defined by \((a\leftarrow
b)\Leftrightarrow (b\rightarrow a)\).  We will often have reason to
investigate properties of a relation and its inverse, and it is
worth noting that the inverse of a reflexive relation is reflexive,
and the inverse of a transitive relation is transitive.

If \(A\) has a (partial) multiplication with domain \(D\), then an
equivalence relation \(\sim\) on \(A\) is a {\itshape congruence
relation} (with respect to the multiplication) if \(aDb\),
\(a'Db'\), \(a\sim a'\) and \(b\sim b'\) always implies
\((ab)\sim(a'b')\).  This is exactly what is needed to induce a
(partial) multiplication on the set \(A/\sim\) of equivalence
classes from the multiplication on \(A\) and give a homomorphism
from \(A\) to \(A/\sim\).  If \(A\) is a group, then an equivalence
relation is a congruence relation if and only if the equivalence
class \(N\) of the identity is a normal subgroup of \(A\) and the
relation is precisely the relation ``belongs to the same coset of
\(N\).''

It is convenient to compare relations under containment.  If \(R\)
and \(S\) are binary relations on a set \(A\), then \(R\subseteq S\)
has the same meaning as for all \((a,b)\in A\times A\), \(aRb\)
implies \(aSb\).

\subsection{Distinguished representatives} We now deal with
properties of relations that lead easily to distinguished
representatives in equivalence classes.  This material is covered
with varying detail and terminology in \cite{squier:word-problems},
\cite{cohen:cgt}, \cite{cohen:rewrite}, \cite{der-jou:rewrite} and
\cite{BKR:reductions}.  Our terminology is closest to
\cite{der-jou:rewrite}.  We assume a binary relation \(\rightarrow\)
on a set \(A\) and we let \(\sim\) be the equivalence relation
generated by \(\rightarrow\).

We let \(\overset{=}\rightarrow\) denote the reflexive closure of
\(\rightarrow\), we let \(\overset{+}{\rightarrow}\) denote the
transitive closure of \(\rightarrow\), and we let
\(\overset{*}\rightarrow\) denote the reflexive, transitive closure
of \(\rightarrow\).  We use \(\overset{=}\leftarrow\),
\(\overset{+}\leftarrow\) and \(\overset{*}\leftarrow\) for the
inverses, respectively, of \(\overset{=}\rightarrow\),
\(\overset{+}\rightarrow\) and \(\overset{*}\rightarrow\).

We use combinations of symbols to represent compositions.  For
example, \(\overset{*}\leftarrow \circ \overset{*}\rightarrow\) is
the relation verbally represented by \(a \overset{*}\leftarrow \circ
\overset{*}\rightarrow b\) if and only if there is a \(c\in A\) so
that \(a\overset{*}\leftarrow c\) and \(c\overset{*}\rightarrow b\).
We can now define the following properties of \(\rightarrow\)
\begin{enumerate} \item {\itshape Church-Rosser}: \(\sim
\,\,\,\subseteq\,\,\, \overset{*}\rightarrow \circ
\overset{*}\leftarrow\). \vspace{3pt} \item {\itshape confluent}:
\(\overset{*}\leftarrow \circ \overset{*}\rightarrow
\,\,\,\subseteq\,\,\, \overset{*}\rightarrow \circ
\overset{*}\leftarrow\).  \vspace{3pt} \item {\itshape strongly
confluent}: \(\leftarrow \circ \rightarrow \,\,\,\subseteq\,\,\,
\overset{=}\rightarrow \circ \overset{=}\leftarrow\).  \vspace{3pt}
\item {\itshape locally confluent}: \(\leftarrow \circ \rightarrow
\,\,\,\subseteq\,\,\, \overset{*}\rightarrow \circ
\overset{*}\leftarrow\).  \end{enumerate} The terminology above is
not applied uniformly in all references.  The implications
\[(3)\Rightarrow (2) \Leftrightarrow (1) \Rightarrow (4)\] are easy
to show.  The implication \((4)\Rightarrow(1)\) is false in general
(see \cite{BKR:reductions} and \cite{der-jou:rewrite}) and a
condition that makes it true is our next interest.

We say \(\rightarrow\) is {\itshape terminating} if there is no
infinite sequence \((a_0, a_1, \ldots)\) with \(a_i\rightarrow
a_{i+1}\) and \(a_i\ne a_{i+1}\) for all \(i\ge0\).  If \(B\subseteq
A\), then a {\itshape minimal} element of \(B\) is some \(b\in B\)
so that \(a\in B-\{b\}\) implies that \(b\rightarrow a\) is false.
We say \(\rightarrow\) is {\itshape well founded} if every non-empty
\(B\subseteq A\) has a minimal element.  A minimal element of \(A\)
is called {\itshape irreducible}.  It is easy that \(\rightarrow\)
is terminating if and only if it is well founded.  A subset \(B\) of
\(A\) is {\itshape inductive} if for each \(a\in A\), we have \(a\in
B\) whenever \(\{x\mid a \rightarrow x,\,\,x\ne a\}\subseteq B\).
The Principle of Noetherian Induction (see \cite{cohn:univ-algebra})
states: {\itshape If \(\rightarrow\) is well founded (equivalently,
terminating), then any inductive subset of \(A\) equals \(A\).}
This is seen by noting that the complement of an inductive subset
has no minimal element.

Noetherian induction gives a not very hard proof of the following.

\begin{prop}[\cite{newman:comb}]\mylabel{NewmanLemma} For a
terminating binary relation \(\rightarrow\) on a set \(A\) and its
generated equivalence relation \(\sim\), the following are
equivalent.  {\claimenum \item \(\rightarrow\) is
Church-Rosser.  \item \(\rightarrow\) is confluent.  \item
\(\rightarrow\) is locally confluent.  \item Each equivalence class
under \(\sim\) contains a unique irreducible element of \(A\).
\end{enumerate}} If (a)--(d) hold, then for \(a\in A\) we have
\(a\overset{*}\rightarrow a'\) where \(a'\) is the unique
irreducible in \(A\) with \(a\sim a'\).  \end{prop}

Proofs, modulo Noetherian induction, are given in
\cite{squier:word-problems} and \cite{BKR:reductions} and are
indicated or given in exercises in \cite{der-jou:rewrite} and
\cite{cohen:cgt}.  

The relation \(\rightarrow\) is called {\itshape complete} if it is
terminating and satisfies any of (a)--(d) from Proposition
\ref{NewmanLemma} above.  The property (d) above of complete
relations will be used repeatedly in this paper.

The practical importance of Proposition \ref{NewmanLemma} is the
equivalence of (c) and (d) since local confluence will often be easy
to verify.  In some cases, strong confluence will hold and the
easier implication ``strong confluence implies confluence'' will
then apply.  However, we will still need terminating to get
irreducible elements and the force of Proposition \ref{NewmanLemma}
will still hold even if its proof can be avoided.

\subsection{Presentations} Using congruence relations, we can
discuss presentations of semigroups and monoids.  We let
\(\N=\{0,1,2,\ldots\}\) be the natural numbers where \(0\) is
regarded as the empty set, \(1=\{0\}\), \(2=\{0,1\}\), and so forth.
We let \(\N'=\N-\{0\}\).  A finite sequence is a function with
domain an element of \(\N\).  The domain \(0\) gives a zero length
sequence (which is equal to the empty function which in turn equals
the empty set).

If \(X\) is a set, then \(X^*\) is the set of finite sequences in
\(X\) with domain from \(\N\), and \(X^+\) is the set of finite
sequences with domain from \(\N'\).  Elements of \(X^*\) or \(X^+\)
are called {\itshape words} in the alphabet \(X\).  We have the
binary operation of concatenation of finite sequences in which \(u\)
and \(v\) with domains \(i\) and \(j\) repectively form \(uv\) with
domain \(i+j\) so that \((uv)_m=u_m\) if \(m<i\) and \(v_{m-i}\)
otherwise.  This gives \(X^*\) an identity element, the unique zero
length sequence, which is absent from \(X^+\).  This identity, being
the empty set, will be denoted \(\phi\) for now, but might be
denoted 1 if appropriate.  The need for care will be discussed
below.

The set \(X^*\) is usually called the {\itshape free monoid} on (the
alphabet) \(X\), and \(X^+\) is usually called the {\itshape free
semigroup} on \(X\).  It is clear that \(X^*\) is a monoid and that
\(X^+\) is a semigroup.  We will not be too concerned with, and will
not review the ``free'' aspect.

A presentation is a pair, usually written \(\langle X\mid R\rangle\)
where \(X\) is a set and \(R\) is a binary relation on either
\(F=X^*\) or \(F=X^+\).  The set \(X\) is called the generating set
and \(R\) is called the relation set.  If the presentation is that
of a monoid, then \(F=X^*\) and if the presentation is that of a
semigroup, then \(F=X^+\).  In both cases, \(F\) has concatenation
as a binary operation.

From \(R\), we build another relation \(R'\) on \(F\) consisting of
all pairs \((xuy, xvy)\) for which all of \(x\), \(u\), \(v\) and
\(y\) are in \(F\) and \((u,v)\) is in \(R\), and then let \(\sim\)
be the equivalence relation generated by \(R'\).  It is a triviality
that \(\sim\) is a congruence relation and we say that the monoid or
semigroup \(F/\sim\) is the one determined by the presentation.

Note that in a monoid presentation, the class of the empty word will
become the global identity of the monoid.  Thus if some non-trivial
word in the generators is to be the identity of the monoid, then the
empty word must appear as one of the components of a pair in the
relation set \(R\).

A presentation of a group is a presentation \(\langle X\mid
R\rangle\) which determines a group by passing through a monoid
presentation (never mentioned) that is derived from \(\langle X\mid
R\rangle\) consisting of \(\langle X \cup X'\mid R \cup I\rangle\)
where \(X'\) is a set disjoint from \(X\) with a fixed one-to-one
correspondence from \(X\) to \(X'\) usually written \(x\mapsto
x^{-1}\) and \(I\) is the set of all pairs of the form
\((xx^{-1},\phi)\) or \((x^{-1}x,\phi)\).  It is a straightforward
exercise to show that this is equivalent to the usual description of
a group presentation.

\subsection{Rewriting rules} Note that if the relation \(R\) in the
above discussion is replaced by its symmetric closure, then the
monoid (semigroup, group) determined by the presentation \(\langle
X\mid R\rangle\) does not change.  Thus one can symmetrize \(R\),
but the non-symmetric version may be useful in its own right.  When
the non-symmetric aspects of a relation set \(R\) are being
emphasized, then the pairs \((u,v)\) in \(R\) are usually called
{\itshape rewriting rules} and are often written as \(u\rightarrow
v\).  Confusingly, the derived relation \(R'\) described above is
then viewed as an extension of the relation \(R\) or \(\rightarrow\)
and is often denoted also by \(\rightarrow\).

We are most interested in completeness.  Let \(\langle X\mid
R\rangle\) be a monoid (semigroup, group) presentation.  If the
derived relation \(R'\) is complete, then we will also say that the
original rewriting rules are complete and also that the presentation
is complete.  From Proposition \ref{NewmanLemma}, if \(\langle X\mid
R\rangle\) is a complete presentation, then each element of the
presented object has a unique representative that is irreducible
(with respect to \(R'\)).  In this situation we say that the
irreducible representative is in {\itshape normal form}.

A rather extreme example of a complete presentation of a semigroup
\(S\) is the presentation \(\langle S\mid M\rangle\) where \(M\)
consists of all the pairs \((s_1s_2\rightarrow s_3)\) where
\(s_1s_2=s_3\) holds in \(S\).  In words, the generating set is all
of \(S\) and the relation set is the multiplication table for \(S\).

The completeness of the presentation is easy.  It is obviously
terminating since applications of \(\rightarrow\) decrease length.
Local confluence breaks into two cases, one of which is trivial and
the other of which reduces immediately to the associativity of the
multiplication.

We need a slight variation for monoids.  The identity of a monoid
presentation is the empty word \(\phi\).  If we were to use the
presentation above for a monoid \(S\) with identity 1, we would have
\(\phi1=1\phi=1\), but \(\phi\) and \(1\) would be distinct.  For
this reason, we let \(S'\) denote \(S-\{1\}\) for a monoid \(S\)
with identity 1.  Now the monoid \(S\) is presented (as a monoid) by
\(\langle S'\mid M\rangle\) where \(M\) is now the multiplication
table as given above with all appearances of \(1\) replaced by
\(\phi\).

A group is a monoid so the above paragraph should apply to groups,
but we do not treat group presentations in the same way as monoid
presentations.  For a group \(G\), the presentation \(\langle G\mid
M\rangle\) with \(M\) the multiplication table as defined for
semigroups gives a valid presentation for \(G\).  This is because
\(\phi\) and the identity 1 for \(G\) are provably in the same
class.

\section{Products of algebraic	
structures\protect\mylabel{PartMulMutual}}

\subsection{Motivation} The structure of semidirect products is well
known.  If a group \(G\) is a product \(AB\) of two subgroups with
\(A\) normal and \(A\cap B=\{1\}\), then conjugation of \(A\) by
elements of \(B\) gives an action of \(B\) on \(A\) by
automorphisms.  If, on the other hand, \(A\) and \(B\) are groups
not known to be subgroups of another group and there is an action of
\(B\) on \(A\) by automorphisms, then a group structure (the
semidirect product) on the set \(A\times B\) can be defined so that
conjugation of \(A\times\{1\}\) by elements of \(\{1\}\times B\)
mirrors the given action.

A similar situation exists even when \(A\) is not assumed to be
normal.  We first look at groups to get the main ingredients.  Let
\(G\) be a group with identity \(1\), and with subgroups \(U\) and
\(A\) satsifying \(U\cap A=\{1\}\) and \(G=UA\).  Then each \(g\in
G\) is uniquely expressible as \(g=u\alpha\) with \(u\in U\) and
\(\alpha\in A\).

We are now in a position to reverse certain products.  With \(u\in
U\) and \(\alpha\in A\), consider \(\alpha u\in G\).  We must have
uniqe elements \(u'\in U\) and \(\alpha'\in A\) so that \(\alpha
u=u'\alpha'\).  This defines two functions \((\alpha,u)\mapsto
\alpha^u\in A\) and \((\alpha,u)\mapsto \alpha\cdot u\in U\) that
are unique in that they satisfy \(\alpha u=(\alpha\cdot
u)(\alpha^u)\) for all \(u\in U\) and \(\alpha\in A\).

\subsection{Generalization} In the discussion above, the assumption
\(U\cap A=\{1\}\) was only used to conclude that each \(g\in G\) was
uniquely a product of the form \(ua\) with \(u\in U\) and \(a\in
A\).  We generalize by making this conclusion a hypothesis.  We also
generalize by removing the requirement that the algebraic structure
be a group.  The following triviality is worth stating separately
since it allows us to make a useful definition.

\begin{lemma} Let \((M,D,m)\) be a set with a partial multiplication
and let \(U\) and \(A\) be subsets of \(M\).  Assume that every
\(x\in M\) is uniquely expressible in the form \(x=u\alpha\) with
\(u\in U\) and \(\alpha\in A\).  Then there are functions \((A\times
U)\cap D\rightarrow A\) written \((\alpha,u)\mapsto \alpha^u\), and
\((A\times U)\cap D\rightarrow U\) written \((\alpha,u)\mapsto
\alpha\cdot u\) defined by the property that \(\alpha u=(\alpha\cdot
u)(\alpha^u)\) whenever \(\alpha Du\).  \end{lemma}

In the setting of the above lemma, the functions \((\alpha,
u)\mapsto \alpha^u\) and \((\alpha, u)\mapsto \alpha\cdot u\)
defined on \((A\times U)\cap D\) will be called the {\itshape mutual
actions defined by the multiplication}.  These actions are
``internally'' generated by the multiplication.  We can also impose
actions ``externally.''

Let \((U, D_U, m_U)\) and \((A, D_A, m_A)\) be sets with partial
multiplications and let \(H\sub A\times U\).  Assume there are
functions \(H\rightarrow A\) written \((\alpha ,u)\mapsto \alpha
\cdot u\) and \(H\rightarrow U\) written \((\alpha ,u)\mapsto \alpha
^u\) defined on \(H\).  We will call such functions {\itshape mutual
actions} defined on H.  However, unless the mutual actions possess
some useful properties, little can be done.  The next lemma
establishes some properties of (useful) mutual actions.  (Recall
that a categorical multiplication is also associative by
definition.)

\begin{lemma}\mylabel{MutualActs} Let \((M,D,m)\) be a set with a
categorical partial multiplication and let \(U\) and \(A\) be
subsets of \(M\) that are closed under the multiplication.  Assume
that every \(x\in M\) is uniquely expressible in the form
\(x=u\alpha\) with \(u\in U\) and \(\alpha\in A\), and let
\((\alpha,u)\mapsto \alpha^u\) and \((\alpha,u)\mapsto \alpha\cdot
u\) defined on \((A\times U)\cap D\) be the mutual actions defined
by the multiplication.  Let \(\alpha\) and \(\beta\) come from \(A\)
and \(u\) and \(v\) come from \(U\).  Then in each of the following,
if one side is defined, then the other side is defined and the two
sides are equal. {\claimenum \item \((\alpha\beta)\cdot u =
\alpha\cdot (\beta\cdot u)\).  \item \((\alpha\beta)^u =
\alpha^{(\beta\cdot u)}\beta^u\).  \item \(\alpha\cdot (uv) =
(\alpha\cdot u)(\alpha^u\cdot v)\).  \item \(\alpha^{(uv)} =
(\alpha^u)^{v}\).  \claimenumend} \end{lemma}

\begin{proof} If the left side of one of  (a) or (b) is defined,
then the left side of the other of (a) or (b) is defined and we have
the equality \((\alpha\beta)u=((\alpha\beta)\cdot u)(\alpha\beta)^u\).

But if the left side of one of (a) or (b) is defined, then \(\alpha
D\beta\) and \((\alpha\beta)Du\) hold.  We have the following chain
of equalities with each step using right or left associativity or a
hypothesis of the lemma as appropriate: \[\begin{split}
(\alpha\beta)u &= \alpha(\beta u) = \alpha((\beta\cdot
u)(\beta^u)) =(\alpha(\beta\cdot u))\beta^u \\ &=((\alpha\cdot
(\beta\cdot u))\alpha^{(\beta\cdot u)})\beta^u = (\alpha\cdot
(\beta\cdot u))(\alpha^{(\beta\cdot u)}\beta^u). \end{split}\]  In
the last expression the left factor lies in \(U\) and the right
factor lies in \(A\).  The equalities in (a) and (b) follow from the
uniqueness hypotheses of the lemma.  

The rest of the argument is left to the reader.  The argument based
on the assumption that the right side of one of (a) or (b) is
defined will need the fact that the multiplication is categorical.
The arguments relating to (c) and (d) are similar to those for (a)
and (b).  \end{proof}

Note that the hypotheses of Lemma \ref{MutualActs} apply
non-trivially to the free abelian monoid on two generators, but not
to the free abelian semigroup on two generators.

Commuting elements act nicely.  The following lemma and corollary
are trivial.

\begin{lemma} Assume the hypotheses and notation of Lemma
\ref{MutualActs} without assuming that the multiplication is
categorical.  Assume that \(u\in U\) and \(\alpha\in A\) commute.
Then \(\alpha\cdot u=u\) and \(\alpha^u=\alpha\).  \end{lemma}

\begin{cor}\mylabel{HowOneActs} Assume the hypotheses and notation
of Lemma \ref{MutualActs} without assuming that the multiplication
is categorical.  Assume that \(M\) has a global identity 1.  If
\(1\in U\), then for all \(\alpha\in A\), we have \(\alpha\cdot
1=1\) and \(\alpha^1=\alpha\).  If \(1\in A\), then for all \(u\in
U\), we have \(1\cdot u=u\) and \(1^u=1\).  \end{cor}

In the situation described by Lemma \ref{MutualActs}, we will refer
to the multiplication on \(M\) as the {\itshape \ZS} product of
\(U\) and \(A\) and write it as \(U\zapprod A\).  Lemma
\ref{MutualActs} only partly justifies this terminology.  Lemma
\ref{ZappaRecon} below finishes the justification.  However, this
will only give an ``internal'' product.  The situation with the
``external'' product is more complicated.  The next lemma will
partially justify some of the assumptions we will make when we
define an external version of the product.

\begin{lemma}\mylabel{ExistZSIdents} Assume the hypotheses and
notation of Lemma \ref{MutualActs}.  If for \(u\in U\) there is a
\(\beta\in A\) for which \(u\beta\) exists, then there is an
\(\alpha_u\in A\) that is a left identity for \(A\) for which
\(u=u\alpha_u\).  If for \(\alpha\in A\) there is a \(v\in U\) for
which \(v\alpha\) exists, then there is a \(u_\alpha\in U\) that is
a right identity for \(U\) for which \(\alpha=u_\alpha \alpha\).
\end{lemma}

\begin{proof} We must have \(u=u_1\alpha_u\) with \(u_1\in U\) and
\(\alpha_u\in A\).  Now \(u\beta = (u_1\alpha_u)\beta =
u_1(\alpha_u\beta)\) and our assumptions imply \(u=u_1\).  Thus
\(u=u_1\alpha_u=u\alpha_u\).  Now consider \(\alpha\in A\) with
\(\alpha_u\alpha\) defined.  Since the multiplication is
categorical, we have that \(u(\alpha_u\alpha)\) is defined and
equals \((u\alpha_u)\alpha = u\alpha\) by associativity.  Again, our
assumptions imply that \(\alpha_u\alpha=\alpha\).  The rest of the
argument is similar to what we have given.  \end{proof}

It turns out that an external product based on (externally imposed)
mutual actions can be defined in extreme generality.  A (poorly
behaved) multiplication emerges even if no additional assumptions
are made.  In general there is not a perfect correspondence between
the internal and external products.  For one thing, embeddings of
the factors might not exist in an external product.  Under extra
hypotheses, the external product can be made to correspond to an
internal product, but this removes considerable useful generality.
Even though the internal product is not fully justified, and an
external product has not yet been defined, we will refer to the \ZS
product as if it has a definition.

\subsection{An example, part I}\protect\mylabel{ZappaExmpl} We give
an example that shows the applicability of the \ZS product to sets
with multiplications that are not full.  It will be revisited in
Section \ref{ZappaExmplII}.

Our binary operations will be compositions on the morphisms of
categories.  Thus we will treat the symbols that represent the
categories as if they primarily represent the collections of
morphisms.

We will start with too much generality and reduce it in stages.

The multiplication on the morphisms of a category that is given by
the composition is a categorical partial
multiplication.\footnote{Not all categorical multiplications arise
in this way.  For a categorical multiplication \((M,D,m)\) to arise
from a category, it must also have full identities and satisfy the
{\itshape digraph rule}: \(aDb\), \(cDb\) and \(cDd\) implies
\(aDd\).}  Thus if we have a category \(G\) and subcategories \(U\)
and \(A\) so that every morphism in \(G\) is uniquely expressible as
\(u\alpha\) with \(u\in U\) and \(\alpha\in A\), then we can write
\(G=U\zapprod A\).  This is little more than a repetition of Lemma
\ref{MutualActs}.

We now reduce the generality the first time.

Let \(G\) be a category and for each object \(x\) in \(G\), let
\(U_x\) be a subgroup of the monoid \(G_x\).  Let \(\widehat{U}\)
denote the totally disconnected groupoid that is the union of all
the \(U_x\).

Further assume that there is a subcategory \(A\) of \(G\) so that
every morphism of \(G\) is uniquely expressible as \(u\alpha\) with
\(u\in \widehat{U}\) and \(\alpha\in A\).  From Lemma
\ref{MutualActs}, we get the mutual actions \((\alpha,u)\mapsto
(\alpha\cdot u)\in \widehat{U}\) and \((\alpha,u)\mapsto
(\alpha^u)\in A\) on a subset of \(A\times \widehat{U}\) satisfying
the identities as given in that lemma and we can write
\(G=\widehat{U}\zapprod A\).

We can say more.  First we note that for all objects \(x\), we have
\(1_x\) in both \(\widehat{U}\) and \(A\).  Thus all of \((1_x)^u\),
\(1_x\cdot u\), \(\alpha^{1_x}\) and \(\alpha\cdot 1_x\) make sense
if \(u\in U_x\) and \(\alpha\) has domain \(x\).  We get the
following equalities: \((1_x)^u = 1_x\), \(1_x\cdot u = u\),
\(\alpha^{1_x} = \alpha\) and \(\alpha \cdot 1_x = 1_y\) where
\(\alpha:x\rightarrow y\).  Thus not all of the actions associated
with the identities of \(G\) are trivial.

For an \(\alpha\) in some \(G_x\), we must have a unique
factorization \(\alpha=u\beta\) with \(\beta\in A_x\) and \(u\in
U_x\).  Thus each \(G_x\) is of the form \(U_x\zapprod A_x\).  Even
if \(G\) is a connected groupoid there can be significant
irregularity.  If \(G\) is a connected groupoid and \(A\) is a
connected subgroupoid, then all the \(G_x\) are isomorphic and all
the \(A_x\) are isomorphic.  The next paragraph shows that not all
the \(U_x\) need be isomorphic.

Let \(G\) be a connected groupoid, let \(A\) be a connected
subgroupoid that contains all the objects of \(G\), and for each
object \(x\) let \(U_x\) be a subgroup of \(G_x\) so that
\(G_x=U_xA_x\) and \(U_x\cap A_x = \{1_x\}\).  Let \(\widehat{U}\)
denote the totally disconnected groupoid that is the union of all
the \(U_x\).  We leave it as an exercise for the reader to show that
each morphism in \(G\) is uniquely expressible as \(u\alpha\) with
\(u\in \widehat{U}\) and \(\alpha\in A\).  Thus it is easy to build
examples by just choosing an appropriate \(U_x\) independently for
each \(x\).  Permutation groups show that this can be done so that
not all the \(U_x\) are isomorphic.  Let \(S_3\) and \(S_4\) be the
permutation groups on, respectively, the sets \(\{0,1,2\}\) and
\(\{0,1,2,3\}\) with \(S_3\) viewed as a subgroup of \(S_4\) in the
usual way.  Then either the Klein four group or the cyclic group of
order four can be embedded in \(S_4\) as a complement to \(S_3\).

If we cut down on the irregularity by assuming that all the \(U_x\)
are isomorphic, there might still be some irregularity remaining.
It is possible to have a group \(H\) with subgroups \(J\), \(K\) and
\(L\) so that \(J\) and \(L\) are isomorphic complements of the
subgroup \(K\), but no automorphism of \(H\) carries \(J\) to \(L\).
This can even be done with \(K\) normal in \(H\).  For example, let
\(H=S_3\times \Z_2\) where \(S_3\) is the permutation group on three
objects, and \(\Z_2\) is the cyclic group of order 2.  If \(x\) is
an involution in \(S_3\) (which is necessarily not central), and
\(y\) is the generator of \(\Z_2\), then we can let \(K=S_3\times
\{1\}\), \(J=\langle (x,y)\rangle\) and \(L=\langle (1,y)\rangle\).
Now \(J\) and \(L\) are of order 2, but are not equivalently
embedded in \(H\) since one is central and the other is not.  They
are both complements of the normal subgroup \(K\).

Thus in the connected groupoid setting, if all the \(U_x\) are
required to be isomorphic, the pairs \((G_x, U_x)\) might not be
ismorphic as pairs even though all the pairs \((G_x, A_x)\) will be
isomorphic as pairs.  In spite of this, nice things happen.  We now
give the final version of our example.

We now return to an arbitrary category \(G\) with subcategory \(A\)
and assume that a fixed group \(U\) is given.  For each object \(x\)
in \(G\), let \(\phi_x\) be a monomorphism from \(U\) into the
monoid \(G_x\) with image \(U_x\).  Let \(\widehat{U}\) denote the
totally disconnected groupoid that is the union of all the \(U_x\)
and assume that every morphism of \(G\) is uniquely expressible as
\(u\alpha\) with \(u\in \widehat{U}\) and \(\alpha\in A\).  As
above, we get mutual actions \((\alpha,u)\mapsto (\alpha\cdot u)\in
\widehat{U}\) and \((\alpha,u)\mapsto (\alpha^u)\in A\) on a subset
of \(A\times \widehat{U}\) and an expression of \(G\) as the \ZS
product \(G=\widehat{U}\zapprod A\).  We also get the facts
mentioned above about the behavior of the mutual actions on the
identities.

We now extract more information.  Since we have \(\alpha u =
(\alpha\cdot u)(\alpha^u)\) and elements of \(\widehat{U}\) have
equal domain and codomain, the domain and codomain of \(\alpha^u\)
are equal to, respectively, the domain and codomain of \(\alpha\).
We can now define mutual actions between \(A\) and \(U\) instead of
between \(A\) and \(\widehat{U}\).  Define \((\alpha, u)\mapsto
\alpha \cdot u \in U\) and \((\alpha, u)\mapsto \alpha ^u\in A\) on
\(A\times U\) as follows.
\mymargin{ConvertZSP}\begin{equation}\label{ConvertZSP}
\mathrm{For\,\,\,} \alpha:x\rightarrow y \mathrm{\,\,and\,\,} u\in
U, \mathrm{\,\,set\,\,\,} \alpha^u = \alpha^{\phi_x(u)}
\mathrm{\,\,and\,\,} \alpha \cdot u = \phi_y^{-1}(\alpha\cdot
(\phi_x(u)))\end{equation} where the functions on \(A\times
\widehat{U}\) are used where appropriate.  These actions are clearly
well defined and they are also defined on all of \(A\times U\).
Note that for \(u\in U\), we also have \(\alpha^u:x\rightarrow y\)
if \(\alpha\) goes from \(x\) to \(y\).

We leave it as an exercise for the reader to show that these mutual
actions on \(A\times U\) satisfy the four identities (a)--(d) from
Lemma \ref{MutualActs} and also satisfy the four equalities
\(1_x\cdot u=u\), \((1_x)^u = 1_x\), \(\alpha \cdot 1=1\) and
\(\alpha^1=\alpha\) as \(u\), \(\alpha\) and \(x\) run over,
respectively, all \(U\), morphisms of \(A\) and objects of \(A\).
We will see in Section \ref{ZappaExmplII} that this will allow us to
also express \(G\) ``externally'' as \(U\zapprod A\).  The
``twisting'' of the actions that allows for lack of normality in the
setting of groups is flexible enough in this situation to allow for
the inequivalent embeddings of the various \(U_x\) in the \(G_x\).

\subsection{One parameter families of functions} We list some
properties that mutual actions might have and introduce terminology
to go with the properties.  To keep notation down, we will never
refer to these functions with a symbol and will either use
\((\alpha,u)\mapsto \alpha^u\) or \((\alpha,u)\mapsto \alpha\cdot
u\) or similar notation to refer to the functions.  We will
concentrate on the fact that these functions can be viewed as
families of fuctions from a set to itself parametrized over a second
set.

A (partial) function \((A\times U)\rightarrow A\) will be refered to
as a (partial) {\itshape one parameter family} of functions
{\itshape over the parameter set} \(U\) and {\itshape base set}
\(A\).  The family is {\itshape partial} if the domain of the
function is actually a proper subset of \(A\times U\).  Otherwise
the family is {\itshape full}.  If the designation full or partial
is omitted, then it will either be clear from the context or
irrelevant.  In the (partial) function \((A\times U) \rightarrow
U\), the set \(A\) is the parameter set, and \(U\) is the base set.

We will make some definitions for \((A\times U)\rightarrow A\)
written as \((\alpha,u)\mapsto \alpha^u\).  The reader can make the
corresponding definitions for \((A\times U)\rightarrow U\).

The family is {\itshape a family of injections} if \(\alpha=\beta\)
whenever \(\alpha^u=\beta^u\).  The family is {\itshape a family of
surjections} if for each \(u\in U\) and \(\beta\in A\), there is an
\(\alpha\in A\) so that \(\beta=\alpha^u\).

The family is {\itshape confluent} if for all \(\alpha\in A\) and
\(u\) and \(v\) in \(U\) with \(\alpha^u\) and \(\alpha^v\) defined,
there are \(p\) and \(q\) in \(U\) so that \((\alpha^u)^p =
(\alpha^v)^q\).

We say that the parametrized family is {\itshape multiplicative} if
\(U\) is a set with a partial multiplication so that \(uv\) and
\(\alpha^{(uv)}\) are defined and \((\alpha^u)^v = \alpha^{(uv)}\)
whenever \((\alpha^u)^v\) is defined.  If a multiplicative family is
over a semigroup, then we can refer to the family as a {\itshape
semigroup action} where the term is to emphasize that \(\alpha^1\)
is not necessarily \(\alpha\) even if \(U\) is in fact a monoid.

The parametrized family is {\itshape trivial} if \(\alpha^u\) is
always equal to \(\alpha\) when defined.  This will usually occur
when the family is an action and in this case we will say that the
action is trivial.

\begin{lemma} A full semigroup action over a semigroup with common
right multiples is confluent.  \end{lemma}

\begin{proof} Let the action be written \((\alpha,u)\mapsto
\alpha^u\).  Let \(up=vq\) be a common right multiple of \(u\) and
\(v\).  Now \((\alpha^u)^p=\alpha^{(up)}=
\alpha^{(vq)}=(\alpha^v)^q\).  \end{proof}

A family is {\itshape coconfluent} if whenever \(\alpha^u=\beta^v\),
there are \(\gamma\), \(p\) and \(q\) so that \(\alpha=\gamma^p\)
and \(\beta=\gamma^q\).  A multiplicative action is {\itshape
strongly coconfluent} if whenever \(\alpha^u=\beta^v\) and \(u\) and
\(v\) have a common left multiple, there are \(\gamma\), \(p\) and
\(q\) so that \(\alpha=\gamma^p\), \(\beta=\gamma^q\) and \(pu=qv\).

\begin{lemma}\mylabel{StrongCoConf} Let \(A\times B\rightarrow A\)
written \((\alpha,u)\mapsto \alpha^u\) be a full, strongly
coconfluent multiplicative action.  Assume that \(B\) is a right
cancellative semigroup, assume that \(\alpha^u=\beta^v\) and assume
that \(u\) and \(v\) have a least common left multiple \(l=au=bv\).
Then there is a \(\gamma\in A\) so that \(\alpha=\gamma^a\) and
\(\beta=\gamma^b\).  \end{lemma}

\begin{proof} A least common left multiple is a common left
multiple, so there are \(\delta\), \(p\) and \(q\) so that \(pu=qv\)
and so that \(\alpha=\delta^p\) and \(\beta=\delta^q\).  Since
\(l=au=bv\) is a least common left multiple, there exists a \(k\) so
that \(pu=qv=kl=kau=kbv\).  By right cancellativity, we have
\(p=ka\) and \(q=kb\).  Let \(\gamma=\delta^k\).  We have
\(\gamma^a=\delta^{(ka)}=\delta^p = \alpha\) and
\(\gamma^b=\delta^{(kb)}=\delta^q=\beta\).  \end{proof}

\subsection{Building multiplications} We now start the work on
reversing Lemma \ref{MutualActs}.  We will do so in pieces since the
calculations are lengthy and we will need the pieces separated in
what we do later.  The reader should keep the example in paragraph
\ref{ZappaExmpl} in mind.

Let \((U, D_U, m_U)\) and \((A, D_A, m_A)\) be two sets with partial
multiplications.  Let \(H\) be a subset of \(A\times U\) and assume
that we have mutual actions \((\alpha ,u)\mapsto \alpha \cdot u\)
and \((\alpha ,u)\mapsto \alpha ^u\) defined on \(H\).  Let \(E\) be
a subset of \(U\times A\).

We now define a partial multiplication \(E\times E\rightarrow
U\times A\)
by\mymargin{TheMult}\begin{equation}\label{TheMult}(u,\alpha
)(v,\beta ) = (u(\alpha \cdot v), \alpha ^v\beta )\end{equation} by
which we mean that the left side is defined if and only if all of
\((\alpha, v)\in H\), \((u, \alpha\cdot v)\in D_U\) and \((\alpha^v,
\beta)\in D_A\) hold.  This partial multiplication is well defined
where it is defined.  However, this multiplication possesses very
few other properties unless assumptions are made about the mutual
actions, the partial multiplications on \(U\) and \(A\) and the sets
\(H\) and \(E\).  It is not even guaranteed that the image of this
mulitplication is contained in \(E\).  Often, but not always, we
will let \(E=U\times A\).  We list some ingredients from which we
can build appropriate assumptions.  The identities are the messiest.

{\begin{enumerate}
\renewcommand{\theenumi}{P\arabic{enumi}}

\item \mylabel{ClosureRules} \begin{enumerate} 

\item \((\alpha\cdot u, \alpha^u)\) is in \(E\) whenever \((\alpha,
u)\) is in \(H\).  

\item \((u, \alpha)\in E\) and \((v, u)\in D_U\) imply \((vu,
\alpha)\in E\).

\item \((u, \alpha)\in E\) and \((\alpha, \beta)\in D_A\) imply
\((u, \alpha\beta)\in E\).

\end{enumerate}

\item\mylabel{AssocRules} \begin{enumerate} 

  \item \((\alpha, u)\mapsto \alpha\cdot u\) is multiplicative:
  \((\alpha\beta)\cdot u = \alpha\cdot (\beta\cdot u)\).

  \item  \((\alpha\beta)^u = \alpha^{\beta\cdot u}\beta^u\).

  \item  \(\alpha\cdot (uv) = (\alpha\cdot u)(\alpha^u\cdot v)\).

  \item \((\alpha,u)\mapsto \alpha^u\) is multiplicative:
      \(\alpha^{uv} = (\alpha^u)^v\).  \end{enumerate}

\item \mylabel{CatRules} \begin{enumerate}

\item \((\alpha, \beta)\in D_A\) and \((\beta, u)\in H\)
imply \((\alpha\beta, u)\in H\).

\item \((u,v)\in D_U\) and \((\alpha, u)\in H\) imply \((\alpha,
uv)\in H\).

\end{enumerate}

\item\mylabel{CancRules} \begin{enumerate}

  \item \((\alpha,u)\mapsto \alpha^u\)  is an injective family.

  \item \((\alpha,u)\mapsto \alpha\cdot u\) is an injective family.
      \end{enumerate}

\item\mylabel{SolnRules} \begin{enumerate}

  \item\((\alpha ,u)\mapsto \alpha \cdot u\) is a surjective family.

  \item\((\alpha,u)\mapsto \alpha^u\) is a surjective family.
  \end{enumerate}

\item\mylabel{FullRule} \(H=A\times U\). 

\item\mylabel{IdentRules} \begin{enumerate}

  \item \(\alpha^u = \alpha\) for every
  \(\alpha\in A\) and right identity \(u\) for \(U\).   

  \item \(\alpha^u = \alpha\) for every
  \(\alpha\in A\) and left identity \(u\) for \(U\).   

  \item \(\alpha\cdot u = u\) for every
  \(u\in U\) and right identity \(\alpha\) for \(A\). 

  \item \(\alpha\cdot u = u\) for every
  \(u\in U\) and left identity \(\alpha\) for \(A\). 

  \item If \(u\) is a right identity for \(U\), then \((\alpha,
  u)\in H\) and \(\alpha\cdot u\) is a right identity for \(U\) for
  all \(\alpha\in A\).

  \item If \(\alpha\) is a left identity for \(A\), then \((\alpha,
  u)\in H\) and \(\alpha^u\) is a left identity for \(A\) for all
  \(u\in U\).

  \item \(u\) is a right identity for \(U\) whenever
  \(\alpha\cdot u\) is a right identity for \(U\).

  \item \(\alpha\) is a left identity for \(A\) whenever
  \(\alpha^u\) is a left identity for \(A\).

  \end{enumerate}

\item\mylabel{CoConfRule} \((\alpha,u)\mapsto \alpha^u\) is strongly
coconfluent.

\end{enumerate}}

Note that assumptions \tref{AssocRules}(a--d) read the same as the
four conclusions of Lemma \ref{MutualActs}.  

We can deal with the closure question.

\begin{lemma}\mylabel{PartMulClosure} Let \((U, D_U, m_U)\) and
\((A, D_A, m_A)\) be sets with partial multiplications and with
mutual actions \((\alpha ,u)\mapsto \alpha \cdot u\) and \((\alpha
,u)\mapsto \alpha ^u\) defined on \(H\subseteq A\times U\).  Put the
partial multiplication \tref{TheMult} on \(E\subseteq U\times A\).
Then the following hold.

{\begin{enumerate}
\renewcommand{\theenumi}{{\itshape \roman{enumi}}}

\item If \tref{ClosureRules}(a--c) hold, then the multiplication is
closed on \(E\).

\item If \(\pi_U\) and \(\pi_A\) are the projections of \(U\times
A\) to \(U\) and \(A\), respectively, then the multiplication
\tref{TheMult} is defined on all of \(E\) if \(H\sub \pi_A(E)\times
\pi_U(E)\).

\keepitem \end{enumerate}}  \end{lemma}

\begin{proof} Consider a pair \((u, \alpha)\), \((v, \beta)\) that
\tref{TheMult} applies to.  Since \((\alpha, v)\in H\),
\tref{ClosureRules}(a) implies that \((\alpha \cdot v, \alpha^v)\)
is in \(E\).  Since \((u, \alpha \cdot v)\in D_U\),
\tref{ClosureRules}(b) implies that \((u(\alpha\cdot v), \alpha^v)\)
is in \(E\).  Since \((\alpha^v, \beta)\in D_A\),
\tref{ClosureRules}(c) implies that \((u(\alpha\cdot v),
\alpha^v\beta)\) is in \(E\) This proves (i).  Item (ii) is trivial.
\end{proof}

From this point, \tref{ClosureRules}(a--c) will be included in all
hypotheses.

Note that the multiplication \tref{TheMult} is defined on all of
\(E\) if \tref{FullRule} holds, and on all of \(U\times A\) if
\(E=U\times A\) and \tref{FullRule} holds.  Note further that
\tref{ClosureRules}(a--c) hold if \(E=U\times A\).

\subsection{Associativity} We will assume that pieces of
\tref{AssocRules}(a) through \tref{AssocRules}(d) hold for the
mutual actions \((\alpha ,u)\mapsto \alpha \cdot u\) and \((\alpha
,u)\mapsto \alpha ^u\).  Our assumptions will take the form of
identities with implications.  For each of the four identities
\tref{AssocRules}(a) through \tref{AssocRules}(d), we will get two
assumptions so that we end up with eight possible assumptions.  The
two assuptions obtained from \tref{AssocRules}(a) will be labelled
\((a\Rightarrow)\) and \((a\Leftarrow)\) and the other six
assumptions will be labelled similarly.

We define these eight assumptions by example and say that
\((a\Rightarrow)\) is to mean that if the left side of
\tref{AssocRules}(a) and its intermediate products are defined, then
the right side of \tref{AssocRules}(a) and all its intermediate
products are defined and the equality holds.  The other seven
assumptions are defined similarly.  We will also say, for example,
that \tref{AssocRules}\((a)\) is satisfied if both
\((a\Rightarrow)\) and \((a\Leftarrow)\) are statisfied.  

We have more notation.  If there is a partial multiplication on
\(U\), then we will write \(\rightarrow_U\) to indicate that the
multiplication on \(U\) is right associative and will write
\(\leftarrow_U\) if the multiplication is left associative.  We now
give the lemma.

\begin{lemma}\mylabel{PartMulAssoc} Let \((U, D_U, m_U)\) and \((A,
D_A, m_A)\) be sets with partial multiplications and with mutual
actions \((\alpha ,u)\mapsto \alpha \cdot u\) and \((\alpha
,u)\mapsto \alpha ^u\) defined on \(H\subseteq A\times U\).  Assume
\tref{ClosureRules}(a--c) and put the partial multiplication
\tref{TheMult} on \(E\subseteq U\times A\).  Then the following
hold.

{\begin{enumerate}\useitem
\renewcommand{\theenumi}{{\itshape \roman{enumi}}}

\item If \(U\) and \(A\) are both right associative and
\((a\Rightarrow)\),
\((b\Rightarrow)\),
\((c\Leftarrow)\) and
\((d\Leftarrow)\) all hold, then the partial
multiplication on \(U\times A\) is right associative.

\item If \(U\) and \(A\) are both left associative and
\((a\Leftarrow)\),
\((b\Leftarrow)\),
\((c\Rightarrow)\) and
\((d\Rightarrow)\) all hold, then the partial
multiplication on \(U\times A\) is left associative.

\item If \(U\) and \(A\) are both associative and
\tref{AssocRules}(a--d) all hold, then
the partial multiplication on \(U\times A\) is associative.
\keepitem \end{enumerate}}  \end{lemma}

\begin{proof}   The proof is contained in the following chain of
equalities.  \begin{xxalignat}{3} 
&  &     &((u,\alpha )(v,\beta ))(w,\gamma ) & &(\ref{TheMult}) \\
&(\ref{TheMult}) & = &(u(\alpha \cdot v), \alpha ^v \beta )(w,\gamma
) & &(\ref{TheMult}) \\
&(\ref{TheMult}) & = &((u(\alpha \cdot v))((\alpha ^v\beta )\cdot
w), (\alpha ^v \beta )^w \gamma ) & &(a\Leftarrow) \\
&(a\Rightarrow) & = &((u(\alpha \cdot v))(\alpha ^v \cdot (\beta
\cdot w)), (\alpha ^v\beta )^w \gamma ) & &(\leftarrow_U) \\
&(\rightarrow_U) & = &(u((\alpha \cdot v)(\alpha ^v \cdot (\beta
\cdot w))), (\alpha ^v\beta )^w\gamma ) & &(c\Rightarrow) \\
&(c\Leftarrow) & = &(u(\alpha \cdot(v(\beta \cdot w))), (\alpha
^v\beta )^w\gamma ) & &(b\Leftarrow) \\
&(b\Rightarrow) & = &(u(\alpha \cdot(v(\beta \cdot w))), ((\alpha
^v)^{\beta \cdot w}\beta ^w)\gamma ) & &(d\Rightarrow) \\
&(d\Leftarrow) & = &(u(\alpha \cdot (v(\beta \cdot w))), (\alpha
^{v(\beta \cdot w)}\beta ^w)\gamma ) & &(\leftarrow_A) \\
&(\rightarrow_A) & = &(u(\alpha \cdot (v(\beta \cdot w))), \alpha
^{v(\beta \cdot w)}(\beta ^w\gamma )) & &(\ref{TheMult}) \\
&(\ref{TheMult}) & = &(u,\alpha )(v(\beta \cdot w), \beta ^w\gamma )
& &(\ref{TheMult}) \\
&(\ref{TheMult}) & = &(u,\alpha )((v,\beta )(w,\gamma ))
\end{xxalignat}
The notations on the left give justifications for the equalities if
they are read from top to bottom, and the notations on the right
give justifications for the equalities if they are read from bottom
to top.  The references to \tref{TheMult} are to the definition of
the multiplication.  \end{proof}

\subsection{Reconstruction} We can now prove that \tref{TheMult}
reproduces multiplications in certain situations.

\begin{lemma}\mylabel{ZappaRecon} Let \((M,D,m)\) be a set with a
categorical, partial multiplication and let \(U\) and \(A\) be
subsets of \(M\) that are closed under the multiplication.  Assume
that every \(x\in M\) is uniquely expressible in the form
\(x=u\alpha\) with \(u\in U\) and \(\alpha\in A\), and let
\((\alpha,u)\mapsto \alpha^u\) and \((\alpha,u)\mapsto \alpha\cdot
u\) defined on \(H=(A\times U)\cap D\) be the mutual actions
defined by the multiplication.  Use these mutual actions and
\tref{TheMult} to build a partial multiplication on
\(E=D\cap(U\times A)\).  Then sending \((u,\alpha)\) in \(E\) to
\(u\alpha\) in \(M\) is an isomorphism of partial multiplications.
\end{lemma}

\begin{proof}  The hypotheses immediately make the function one to
one and onto.

Consider two pairs \((u, \alpha)\) and \((v, \beta)\) in \(U\times
A\) and assume that \(\alpha \cdot v\) and \(\alpha^v\) exist.  If
all six of \(uD\alpha\), \(\alpha Dv\), \(vD\beta\), \(u D (\alpha
\cdot v)\), \((\alpha \cdot v)D(\alpha ^v)\) and \((\alpha
^v)D\beta\) hold, then the properities implied by ``categorical''
give \((u\alpha)D(v\beta)\), \((u(\alpha \cdot
v))D((\alpha^v)\beta)\) and
\mymargin{ReconEquality}\begin{equation}\label{ReconEquality}
(u\alpha)(v\beta) = (u(\alpha\cdot
v))((\alpha^v)\beta).\end{equation}

Let \((u, \alpha)\) and \((v, \beta)\) be from \(E\).  From the
conditions built into \tref{TheMult}, the existence of \((u,
\alpha)(v, \beta)\) and the categorical assumption lead to all the
conditions that imply \tref{ReconEquality}.  If
\((u\alpha)D(v\beta)\), then the assumptions behind associativity
also lead to the conditions that imply \tref{ReconEquality}.  Thus
the domains of definition of the multiplications on \(M\) and \(E\)
are carried onto one another, and the product itself is preserved.
\end{proof}

We can now make a definition.  Assuming the setting, hypotheses and
notation of Lemma \ref{ZappaRecon}, we say that \((M, D, m)\) is the
{\itshape (internal) \ZS product} of \(U\) and \(A\) and write
\(M=U\zapprod A\).

Results are not attractive when one attempts to create an
``external'' version that corresponds to the internal version.  The
next lemma is the best that we have been able to do in this
direction.

We take hints from Lemmas \ref{MutualActs}, \ref{ExistZSIdents} and
\ref{ZappaRecon} and try to find conditions that imply that the
product is categorical and that the factors are embedded.  The next
lemma takes care of these as best we can.  First we need two
definitions.

Let \((U, D_U, m_U)\) and \((A, D_A, m_A)\) be sets with partial
multiplications with mutual actions \((\alpha ,u)\mapsto \alpha
\cdot u\) and \((\alpha ,u)\mapsto \alpha ^u\) defined on
\(H\subseteq A\times U\).  Let \(E\sub U\times A\) be given.  A pair
of functions \(i:A\rightarrow U\) and \(j:U\rightarrow A\) will be
called {\itshape embedding functions} if the following properties
hold: \[\begin{gathered} \begin{gathered} j(U)\sub
\{\mathrm{full\,\,identities\,\,of\,\,}A\}, \\ uD_Uv \Rightarrow
j(u)=j(v)=j(uv), \\ (u, j(u))\in E, \\ (j(u), u) \in H, \\ j(u)
\cdot u = u, \\ (j(u))^u = j(u), \\ j(u)D_A j(u),
\end{gathered}\qquad\quad \begin{gathered} i(A)\sub
\{\mathrm{full\,\,identities\,\,of\,\,}U\}, \\ \alpha D_A\beta
\Rightarrow i(\alpha) = i(\beta) = i(\alpha\beta), \\ (i(\alpha),
\alpha) \in E, \\ (\alpha, i(\alpha)) \in H, \\ \alpha\cdot
i(\alpha) = i(\alpha), \\ \alpha^{i(\alpha)} = \alpha, \\ i(\alpha)
D_U i(\alpha), \end{gathered} \\ (\alpha, u)\in H \Rightarrow \big[
i(\alpha)D_U (\alpha \cdot u) \mathrm{\,\,and\,\,} \alpha^u D_A j(u)
\big], \\ \mathrm{and} \\ \begin{aligned} (u, \alpha)\in E
\Rightarrow \big[&(j(u), i(\alpha)) \in H,\,\,j(u)\cdot i(\alpha) =
i(\alpha),\,\, (j(u))^{i(\alpha)} = j(u), \\ &uD_U i(\alpha)
\mathrm{\,\,and\,\,} j(u)D_A \alpha\big]. \end{aligned}
\end{gathered}\]

The reader should note how much simpler the above assumptions get
when \(A\) and \(U\) are monoids, \(E=U\times A\) and \(H=A\times
U\).  This will be formalized below in Lemma
\ref{PartMulPropsII}(\ref{MonoidItem}).

Now let \(f:U\rightarrow U'\) and \(g:A\rightarrow A'\) be functions
(we do not need homomorphisms for this) between sets with
partial multiplication so that there are mutual actions between
\(U\) and \(A\) defined on \(H\sub A\times U\) and also mutual
actions between \(U'\) and \(A'\) defined on \(H'\sub A'\times U'\)
(with the usual notations).  We say that \(f\) and \(g\) {\itshape
preserve the actions} if for all \((\alpha, u)\in H\) we have
\((g(\alpha), f(u))\in H'\), \(g(\alpha^u) = (g(\alpha))^{f(u)}\)
and \(f(\alpha\cdot u) = g(\alpha) \cdot f(u)\).

\begin{lemma}\mylabel{PartMulProd} Let \((U, D_U, m_U)\) and \((A,
D_A, m_A)\) be sets with partial multiplications and with mutual
actions \((\alpha ,u)\mapsto \alpha \cdot u\) and \((\alpha
,u)\mapsto \alpha ^u\) defined on \(H\subseteq A\times U\).  Assume
\tref{ClosureRules}(a--c) and put the partial multiplication
\tref{TheMult} on \(E\subseteq U\times A\).  Then the following
hold.

{\begin{enumerate}\useitem
\renewcommand{\theenumi}{{\itshape \roman{enumi}}}

\item\mylabel{InjItem} If \(i:A\rightarrow U\) and \(j:U\rightarrow
A\) are embedding functions, then sending \(\alpha\in A\) to
\((i(\alpha), \alpha)\) and sending \(u\in U\) to \((u, j(u))\)
gives homomorphic embeddings from \(A\) and \(U\) into the partial
multiplication that \tref{TheMult} puts on \(E\).  Further if the
images of \(A\) and \(U\) under these embeddings are \(A'\) and
\(U'\), respectively, then each element of \(E\) is uniquely
expressible as a product of an element of \(U'\) followed by an
element of \(A'\) and if the mutual actions defined by the
multiplication on \(E\) are placed on \(A'\) and \(U'\), then the
homomorphic embeddings preserve the actions.

\item\mylabel{CatItem} If \(U\) and \(A\) are both categorical and
\tref{AssocRules}(a--d) and \tref{CatRules}(a--b) all hold, then
\tref{TheMult} makes the partial multiplication on \(E\)
categorical.

\keepitem \end{enumerate}}  \end{lemma}

\begin{proof} We first work on \tref{InjItem}.  Sending \(\alpha\in
A\) to \((i(\alpha), \alpha)\) and sending \(u\in U\) to \((u,
j(u))\) clearly give embeddings.  It follows easily from the (heavy
handed) hypotheses that these embeddings are homomorphisms on the
domains of the multiplications.  If \((u, \alpha)\) is in \(E\),
then the equality \((u, j(u))(i(\alpha), \alpha)=(u, \alpha)\) is
equally easy to verify.  The well definedness of the product implies
the uniqueness of the expansion of \((u, \alpha)\) as \((u,
j(u))(i(\alpha), \alpha)\).  Lastly, if \((\alpha, u)\in H\), then a
simple calculation shows that \((i(\alpha), \alpha)(u, j(u)) =
(\alpha \cdot u, \alpha^u)\).  This shows that the embeddings
preserve the mutual actions.

We now work on \tref{CatItem}. From \tref{AssocRules}(a--d) and the
associativity of \(U\) and \(A\), we know that the partial
multiplication on \(E\) is associative.  If both products \((u,
\alpha)(v, \beta)\) and \((v, \beta)(w, \gamma)\) exist with \((u,
\alpha)\), \((v, \beta)\) and \((w, \gamma)\) all in \(E\), then we
must show that \([(u, \alpha)(v, \beta)](w, \gamma)\) exists.  The
existence of \((u, \alpha)[(v, \beta)(w, \gamma)]\) will follow from
an identical argument.

To show that \([(u, \alpha)(v, \beta)](w, \gamma)\) exists, we must
show all three of \((\alpha^v\beta, w)\in H\), \((u(\alpha\cdot v),
(\alpha^v\beta)\cdot w)\in D_U\) and \(((\alpha^v\beta)^w,
\gamma)\in D_A\).  From the existence of \((u, \alpha)(v, \beta)\)
and \((v, \beta)(w, \gamma)\), we have that both \((\alpha, v)\) and
\((\beta, w)\) are in \(H\), both \((u, \alpha \cdot v)\) and \((v,
\beta\cdot w)\) are in \(D_U\) and both \((\alpha^v, \beta)\) and
\((\beta^w, \gamma)\) are in \(D_A\).

Combining \((\alpha, v)\in H\) with \((v, \beta\cdot w)\in D_U\)
with \tref{CatRules}(b) gives \((\alpha, v(\beta\cdot w))\in H\).
This says that \(\alpha\cdot(v(\beta\cdot w))\) exists which implies
by \tref{AssocRules}(c \(\Rightarrow\)) that \(\alpha^v
\cdot(\beta\cdot w)\) exists and \((\alpha^v, \beta\cdot w)\) is in
\(H\).  Combining \((\beta, w)\in H\) and \((\alpha^v, \beta\cdot
w)\in H\) with \tref{AssocRules}(a \(\Leftarrow\)) gives
\((\alpha^v\beta, w)\in H\) which is the first item we must show.

From \tref{AssocRules}(c \(\Rightarrow\)) and the existence of
\(\alpha\cdot(v(\beta\cdot w))\) we also get \((\alpha\cdot v,
\alpha^v\cdot(\beta\cdot w))\in D_U\).  Combining this last fact and
\((u, \alpha\cdot v)\in D_U\) with \tref{AssocRules}(a
\(\Leftarrow\)) and the assumption that \(U\) is categorical gives
\((u(\alpha\cdot v), (\alpha^v\beta)\cdot w)\in D_U\) which is the
second item we must show.

Combining \((\alpha^v, \beta)\in D_A\) and \((\beta, w)\in H\) with
\tref{CatRules}(a) gives \((\alpha^v\beta, w)\in H\) which in turn
implies that \((\alpha^v\beta)^w\) exists.  By \tref{AssocRules}(b
\(\Rightarrow\)), \((\alpha^{v(\beta\cdot w)}, \beta^w)\) is in
\(D_A\).  Combining this last and \((\beta^w, \gamma)\in D_A\) with
the assumption that \(A\) is categorical gives
\((\alpha^{v(\beta\cdot w)}\beta^w, \gamma) \in D_A\).  This is
equivalent to the last item that we must show by \tref{AssocRules}(b
\(\Leftarrow\)).  \end{proof}

The setting and hypotheses needed in both parts of Lemma
\ref{PartMulProd} define the {\itshape external \ZS product} of
\(U\) and \(A\), also written \(U\zapprod A\).  It is unsatisfactory
in the number of hypotheses needed.

There is no easy ``correspondence'' between the internal and
external \ZS products.  We revisit the example of Section
\ref{ZappaExmpl} to establish such a correspondence in a particular
setting.

\subsection{The example, part II}\mylabel{ZappaExmplII} Recall the
example at the end of Section \ref{ZappaExmpl}.  We are given an
arbitrary category \(G\), a subcategory \(A\), a fixed group \(U\),
and for each object \(x\) a monomorphism \(\phi_x\) from \(U\) into
the monoid \(G_x\) with image \(U_x\).  We let \(\widehat{U}\) be
the union of all the \(U_x\) and we assume that every morphism of
\(G\) is uniquely expressible as \(u\alpha\) with \(u\in
\widehat{U}\) and \(\alpha\in A\).

We now take note of the mutual actions \((\alpha, u)\mapsto
\alpha\cdot u\) and \((\alpha, u)\mapsto \alpha^u\) that are defined
by \tref{ConvertZSP} in Section \ref{ZappaExmpl} on all of
\(H=A\times U\).  These mutual actions are defined from the original
actions on \(A\) and \(\widehat{U}\) that are induced from the
multipication on \(G\).  The actions on \(H\) satisfy
\begin{enumerate} \item the four properties \ref{AssocRules}(a--d)
of Lemma \ref{MutualActs}, and \item the equalities \(1_x\cdot
u=u\), \((1_x)^u = 1_x\), \(\alpha \cdot 1=1\) and
\(\alpha^1=\alpha\), \(S(\alpha^u)=S(\alpha)\) and
\(T(\alpha^u)=T(\alpha)\) as \(u\), \(\alpha\) and \(x\) run over,
respectively, \(U\), morphisms of \(A\) and objects of \(A\).
\end{enumerate}

We now let \(E=U\times A\), we fix an object \(x\) in \(G\), and we
define \(i(\alpha)=1\in U\) and \(j(u)=1_x\in A\).  One can verify
that all of the hypotheses of Lemma \ref{PartMulProd} are satisfied.
Thus we are justified in writing \(U\zapprod A\) for the
multiplication that \tref{TheMult} puts on \(U\times A\).

Lemma \ref{ZappaRecon} allows us to conclude that (the
multiplicative structure on the morphisms of) \(G\) has the
structure \(\widehat{U}\zapprod A\).  We wish to show an isomorphism
between \(U\zapprod A\) and \(\widehat{U}\zapprod A\).

Let \(u\) be in \(U\) and \(\alpha:x\rightarrow y\) be in \(A\).  We
send \((u, \alpha)\) in \(U\zapprod A\) to \((\phi_y(u), \alpha)\)
in \(\widehat{U}\zapprod A\).  In \(U\zapprod A\), we have \((u,
\alpha)(v, \beta) = (u(\alpha \cdot v), \alpha^v\beta)\) from
\tref{TheMult}.  Since the domain and codomain of \(\alpha^v\) are
that of \(\alpha\), since \(U\) is a group and since \(H=A\times
U\), we have that \((u, \alpha)(v, \beta)\) exists if and only if
\(\alpha\beta\) exists in \(A\).  Similarly, the images of \((u,
\alpha)\) and \((v, \beta)\) in \(\widehat{U}\zapprod A\) have their
product exist under exactly the same condition.  Thus the mapping
preserves the domain of the product.

From the structure of \(G\) and the fact that sending \((\phi_y(u),
\alpha)\) in \(\widehat{U}\zapprod A\) to \(\phi_y(u)\alpha\) in
\(G\) forms the isomorphism between \(\widehat{U}\zapprod A\) and
\(G\), we get that our mapping is onto.  Our mapping is seen one to
one by considering the domain and ranges of the images of elements,
by noting that every morphism in \(G\) is uniquely an allowable
composition from \(\widehat{U}A\), and from the injectivity of the
various \(\phi_y\).

To see that the mapping preserves products let \(\alpha:x\rightarrow
y\) and \(\beta:z\rightarrow x\), and note that \((u, \alpha)\) and
\((v, \beta)\) in \(U\zapprod A\) are sent to \((\phi_y(u),
\alpha)\) and \((\phi_x(v), \beta)\), respectively, in
\(\widehat{U}\zapprod A\).  In \(U\zapprod A\), the product \((u,
\alpha)(v, \beta)\) is evaluted by \tref{TheMult} as
\((u(\alpha\cdot v), \alpha^v\beta)\) which by definition equals
\((u\phi_y^{-1}(\alpha\cdot \phi_x(v)), \alpha^{\phi_x(v)}\beta)\).
This is sent by our mapping to \[\begin{split}
&(\phi_y[u\phi_y^{-1}(\alpha \cdot \phi_x(v))] ,
\alpha^{\phi_x(v)}\beta) \\ = &(\phi_y(u)(\alpha\cdot \phi_x(v)),
\alpha^{\phi_x(v)}\beta)\end{split}\] which is the product
\((\phi_y(u), \alpha)(\phi_x(v), \beta)\) in \(\widehat{U}\zapprod
A\).

On the other hand, assume we are given a category \(A\) and group
\(U\) with mutual actions between them defined on \(H=A\times U\),
written as above and satisfying (1) and (2) above.  Then the
embedding functions used above allow us to apply Lemma
\ref{PartMulProd} and define \(U\zapprod A\) on \(E=U\times A\).  We
try to make a category from \(U\zapprod A\) by using the objects of
\(A\) as objects of \(U\zapprod A\) and setting
\(S(u,\alpha)=S(\alpha)\) and \(T(u, \alpha)=T(\alpha)\).  Since
\(H=U\times A\), the formula \((u, \alpha)(v, \beta) = (u(\alpha
\cdot v), \alpha^v\beta)\) says that \((u, \alpha)(v, \beta)\)
exists if and only if \(u(\alpha\cdot v)\) exists and
\(\alpha^v\beta\) exists.  Since \(U\) is a group, \(u(\alpha \cdot
v)\) always exists.  Now \(\alpha^v\beta\) exists if and only if
\(\alpha\beta\) exists since \(S(\alpha^v)=S(\alpha)\) and
\(T(\alpha^v)=T(\alpha)\).  Thus \((u, \alpha)(v, \beta)\) exists if
and only if \(\alpha \beta\) exists, and elements of \(U\zapprod A\)
multiply if and only if domain and range match correctly.  An easy
check shows that \((1,1_x)\) acts as an identity at \(x\). Our
attempt to make \(U\zapprod A\) a category is successful.

Now in \(U\zapprod A\), the subcategory
\(\widehat{A}=\{(1,\alpha)\mid \alpha\in A\}\) and subgroupoid
\(\widehat{U}=\{(u,1_x)\mid u\in U, \,\, x\mathrm{\,\,an
\,\,object\,\, in\,\,}A\}\) fit the description of the example of
Section \ref{ZappaExmpl}.  The argument above shows that the
external product \(U\zapprod A\) and the internal product
\(\widehat{U}\zapprod \widehat{A}\) are isomorphic.  From this
isomorphism and Lemma \ref{PartMulProd}(\ref{InjItem}), the mutual
actions between \(U\) and \(A\) deduced from those between
\(\widehat{U}\) and \(\widehat{A}\) as in Section \ref{ZappaExmpl}
are the original assumed actions between \(U\) and \(A\).  We have
essentially shown the following.

\begin{lemma}\mylabel{ZSIntExt} Consider the situations described in
(I) and (II):

(I) A category \(G\) and two subcategories \(A\) and
\(\widehat{U}\) are given where each contains all the objects of
\(G\) and with \(\widehat{U}\) a totally disconnected groupoid so
that every morphism in \(G\) is uniquely expressible as \(u\alpha\)
with \(u\in \widehat{U}\) and \(\alpha\in A\).  Also a group \(U\)
is given with, for each object \(x\) of \(G\), an isomorphism
\(\phi_x:U\rightarrow \widehat{U}_x\).

(II) A category \(A\) and group \(U\) are given with mutual actions
\((\alpha, u)\mapsto \alpha\cdot u\) and \((\alpha, u)\mapsto
\alpha^u\) between them defined on \(H=A\times U\), satisfying
\begin{enumerate} \item the four properties \ref{AssocRules}(a--d)
of Lemma \ref{MutualActs}, and \item the equalities \(1_x\cdot
u=u\), \((1_x)^u = 1_x\), \(\alpha \cdot 1=1\) and
\(\alpha^1=\alpha\), \(S(\alpha^u)=S(\alpha)\) and
\(T(\alpha^u)=T(\alpha)\) as \(u\), \(\alpha\) and \(x\) run over,
respectively, \(U\), morphisms of \(A\) and objects of \(A\)
\end{enumerate} to give a multiplicative structure \(U\zapprod A\)
on all of \(U\times A\).

Given Situation (I), the mutual between \(U\) and \(A\) action
defined by \tref{ConvertZSP} in Section \ref{ZappaExmpl} defines
multiplication on all of \(U\times A\) that fits the description of
Situation (II) so that \(G\) is isomorphic to \(U\zapprod A\).
Given Situation (II), the \ZS product \(U\zapprod A\) is a category
fitting Situation (I).  The processes for passing between the two
situations are inverse to each other.  \end{lemma}

\subsection{Cancellativity and multiples} The next lemma states
conclusions that derive properties that have the word ``right'' in
them.  We leave it to the reader to supply corresponding statements
that have the word ``left.''  We also leave the proofs to the
reader.

\begin{lemma}\mylabel{PartMulProps} Let \((U, D_U, m_U)\) and \((A,
D_A, m_A)\) be sets with partial multiplications and with mutual
actions \((\alpha ,u)\mapsto \alpha \cdot u\) and \((\alpha
,u)\mapsto \alpha ^u\) defined on \(H\subseteq A\times U\).  Assume
\tref{ClosureRules}(a--c) and put the partial multiplication
\tref{TheMult} on \(E\subseteq U\times A\).  Then the following
hold.

{\begin{enumerate}\useitem
\renewcommand{\theenumi}{{\itshape \roman{enumi}}}

\item\mylabel{CancItem} If \(U\) and \(A\) are both right
cancellative and \tref{CancRules}(a) holds, then the partial
multiplication on \(E\) is right cancellative.

\item\mylabel{CRMItem} If \(U\) and \(A\) both have common right
multiples and \tref{FullRule} and \tref{SolnRules}(a) both hold,
then the partial multiplication on \(E\) has common right
multiples.

\item\mylabel{SemigroupItem} If \(A\) and \(U\) are semigroups and
\tref{AssocRules}(a--d) and \tref{FullRule} hold, then
\tref{TheMult} makes \(E\) a semigroup.  \keepitem \end{enumerate}}
\end{lemma}

\subsection{Identities and inverses} Recall that there are many
types of identites possible and even more types of inverses.  See
the definitions before Lemmas \ref{UniqueId} and \ref{FullId}.  What
is given in the next lemma is not hard to prove.  It gives examples
of what can be proven.  From Lemma \ref{UniqueId} we know that some
simplification occurs when associativity is present, and from Lemma
\ref{FullId} we know that the left/right distinction for identities
disappears in a set with a partial multiplication if there are full
identities present.

Recall that a set \(M\) with partial multiplication has left
inverses with respect to right identities if there is an \(x\)
making \(xa=b\) whenever \(b\) is a right identity for \(M\) and
\(ab\) is defined.

In the following, the equality \(E=U\times A\) is deliberate.

\begin{lemma}\mylabel{PartMulPropsII} Let \((U, D_U, m_U)\) and
\((A, D_A, m_A)\) be sets with partial multiplications and with
mutual actions \((\alpha ,u)\mapsto \alpha \cdot u\) and \((\alpha
,u)\mapsto \alpha ^u\) defined on \(H\subseteq A\times U\).  Assume
\tref{ClosureRules}(a--c) and put the partial multiplication
\tref{TheMult} on \(E = U\times A\).  Then the following hold.

{\begin{enumerate}\useitem
\renewcommand{\theenumi}{{\itshape \roman{enumi}}}

\item\mylabel{IdentBaseItem} If and \tref{IdentRules}(a,e,g) hold
and \((u,\alpha)\) is a right identity for \(U\times A\), then
\(u\) is a right identity for \(U\) and \(\alpha\) is a right
identity for \(A\).

\item\mylabel{RIdentItem} If \tref{IdentRules}(a,e) hold then
\((u,\alpha)\) is a right identity for \(U\times A\) whenever
\(u\) is a right identity for \(U\) and \(\alpha\) is a right
identity for \(A\).

\item\mylabel{RIdentsItem} If \(A\) and \(U\) have right identities
and \tref{IdentRules}(a,e,g) and \tref{SolnRules}(a) all hold, then
\tref{TheMult} gives \(U\times A\) right identities.

\item\mylabel{LInvItem} If \(A\) and \(U\) have left inverses with
respect to right identities and \tref{AssocRules}(a,c),
\tref{IdentRules}(a,c,e,g) and \tref{SolnRules}(b) hold, then
\tref{TheMult} gives \(U\times A\) left inverses with respect to
right identities.

\item\mylabel{MonoidItem} If \(A\) and \(U\) are monoids, if
\(H=A\times U\) and if all of \tref{AssocRules}(a--d),
\tref{FullRule} and \tref{IdentRules}(a,d,e,f) all hold, then
\tref{TheMult} makes \(U\times A\) a monoid.  Further, \(j(u)=1_A\)
and \(i(\alpha)=1_U\) make valid embedding functions.

\item\mylabel{GroupItem} If \(A\) and \(U\) are groups, if
\(H=A\times U\) and if all of \tref{AssocRules}(a--d),
\tref{FullRule} and \tref{IdentRules}(a--h) hold, then
\tref{TheMult} makes \(U\times A\) a group.

\keepitem\end{enumerate}}  \end{lemma}

\begin{proof} For \tref{IdentBaseItem}, we have
\((v,\beta)(u,\alpha)=(v,\beta)\) for all \((v,\beta)\) where
\((v,\beta)(u,\alpha)\) is defined.  This says \((v(\beta \cdot u),
\beta^u\alpha)=(v,\beta)\) or \(v(\beta\cdot u)=v\) and
\(\beta^u\alpha=\beta\).  If \(w(\beta\cdot u)\) is defined, then so
is \((w,\beta)(u,\alpha)\) which must equal \((w,\beta)\) and
\(w(\beta\cdot u)=w\).  Thus \(\beta \cdot u\) is a right identity
for \(U\) and from \tref{IdentRules}(g) we have that \(u\) is also a
right identity for \(U\).

Now suppose that \(\gamma\alpha\) is defined.  By
\tref{IdentRules}(e), \(\gamma\cdot u\) is a right identity for
\(U\) and the definition of right identity for \(U\) says that there
is a \(w\in U\) with \(w(\gamma\cdot u)\) defined (and thus equal to
\(w\)).  By \tref{IdentRules}(a), \(\gamma^u=\gamma\) so
\((w,\gamma)(u,\alpha)= (w(\gamma\cdot u), \gamma^u\alpha)= (w,
\gamma\alpha)\) is defined and equals \((w,\gamma)\).  Thus
\(\alpha\) is a right identity for \(A\).

For \tref{RIdentItem}, if \(u\) is a right identity for \(U\) and
\(\alpha\) is a right identity for \(A\), then \(\beta\cdot u\) is a
right identity for \(U\) by \tref{IdentRules}(e) and
\(\beta^u=\beta\) by \tref{IdentRules}(a).  Now for any
\((v,\beta)\) in \(U\times A\) with \((v,\beta)(u,\alpha)\) defined,
we have \[(v,\beta)(u,\alpha)=(v(\beta\cdot u), \beta^u\alpha) =
(v,\beta\alpha)=(v,\beta).\] That there is such a \((v,\beta)\)
follows because there must be a \(\beta\) with \(\beta\alpha\)
defined, and there must be a \(v\) with \(v(\beta\cdot u)\) defined.

In \tref{RIdentsItem}, we are given \((v,\beta)\) and want
\((u,\alpha)\) so that \((v,\beta)(u,\alpha)=(v,\beta)\).  Thus we
want \((v(\beta\cdot u),\beta^u\alpha)=(v,\beta)\).  Since there are
right identities in \(U\), there is a value for \(\beta\cdot u\)
that is a right identity for \(U\) for which \(v(\beta\cdot u)\) is
defined (and thus equal to \(v\)).  From \tref{SolnRules}(a), there
is a \(u\) so that \(\beta\cdot u\) has the right value, and from
\tref{IdentRules}(g), we have that \(u\) is a right identity for
\(U\).  Now \tref{IdentRules}(a) gives that \(\beta^u=\beta\) and
the existence of right identities in \(A\) gives a right identity
\(\alpha\) for \(A\) so that \(\beta^u\alpha=\beta\alpha=\beta\).
Now if \((w,\gamma)(u,\alpha)\) is defined, then we have
\((w(\gamma\cdot u),\gamma^u\alpha) =(w,\gamma)\) since
\tref{IdentRules}(e) gives that \(\gamma\cdot u\) is a right
identity for \(U\) and \tref{IdentRules}(a) gives that
\(\gamma^u=\gamma\) and we already know that \(\alpha\) is a right
identity for \(A\).

In \tref{LInvItem}, we want a solution \((w,\gamma)\) to
\((w,\gamma)(v,\beta) = (u,\alpha)\) when \((u,\alpha)\) is a right
identity for \(U\times A\) and \((v,\beta)(u,\alpha) = (v(\beta\cdot
u), \beta^u\alpha)\) is defined.  Thus we want \((w(\gamma\cdot
v),\gamma^v\beta)=(u,\alpha)\).  From item \tref{IdentBaseItem}, we
know that \(u\) is a right identity for \(U\) and \(\alpha\) is a
right identity for \(A\).  Thus \(\beta^u=\beta\) by
\tref{IdentRules}(a) and \(\beta\alpha\) is defined.  This gives a
solution \(\delta\) to \(\delta\beta=\alpha\).  There is a solution
\(\gamma\) to \(\gamma^v=\delta\) by \tref{SolnRules}(b).  Thus we
have \(\gamma^v\beta=\alpha\).

We want a solution \(w\) to \(w(\gamma\cdot v)=u\).  We know that
\(u\) is a right identity for \(U\), and we need to know that
\((\gamma\cdot v)u\) is defined.  From
\((v,\beta)(u,\alpha)=(v,\beta)\), we know that \(v(\beta\cdot
u)=v\).  Thus \[\gamma\cdot v = \gamma\cdot(v(\beta\cdot u)) =
(\gamma \cdot v)((\gamma^v\beta)\cdot u) = (\gamma\cdot
v)(\alpha\cdot u) = (\gamma\cdot v)u\] where the second equality
uses \tref{AssocRules}(a,c) and the last equality uses
\tref{IdentRules}(c).  Thus \((\gamma\cdot v)u\) is defined (and of
course equal to \(\gamma\cdot v\)), and our hypotheses say that
there is a solution \(w\) to \(w(\gamma\cdot v)=u\).

In \tref{MonoidItem}, \tref{AssocRules} and \tref{FullRule} give an
associative, full multiplication.  The assumed parts of
\tref{IdentRules}, item \tref{RIdentItem} and the unstated variation
of \tref{RIdentItem} for left identities give a full identity.
Since the multiplication is full, we get a global identity.  The
checks that \(i\) and \(j\) are embedding functions is trivial.

In \tref{GroupItem}, we get \tref{SolnRules}(a,b) from the fact that
we are dealing with groups and the rest of the hypotheses say that
\((\alpha,u)\mapsto \alpha\cdot u\) and \((\alpha, u)\mapsto
\alpha^u\) are actions.  The rest follows from \tref{MonoidItem},
\tref{LInvItem}, and the variant of \tref{LInvItem} with left and
right reversed.  \end{proof}

\subsection{Monoids and groups\protect\mylabel{MonGpPars}} Recall
that Lemma \ref{PartMulProd} assembles a long list of hypotheses to
create embeddings of \(A\) and \(U\) into the product created by
\tref{TheMult}.  Global identities tend to make the task much
easier.  If \(A\) has a global identity \(1\), if \(U\times\{1\}
\sub E\), and if \ref{IdentRules}(d,f) hold, then mapping \(u\mapsto
(u,1)\) is a monomorphism and a symmetric statement holds for an
embedding of \(A\).  Under the hypotheses of Lemma
\ref{PartMulPropsII}\tref{MonoidItem} with 1 denoting the global
identities of both \(U\) and \(A\), we have
\((u,1)(1,\alpha)=(u,\alpha)\) and \((1,\alpha)(u,1)= (\alpha\cdot
u, \alpha^u)\).  We can identify \(U\) and \(A\) with
\(U\times\{1\}\) and \(\{1\}\times A\), respectively, and we see
that \(u\alpha\) naturally represents the product in \(U\times A\)
and that the functions \((\alpha,u)\mapsto \alpha\cdot u\) and
\((\alpha,u)\mapsto \alpha^u\) tell how to evaluate \(\alpha u\) by
``reversing'' the product to something of the form \(u'\alpha'\)
with \(u'\in U\) and \(\alpha'\in A\).

One typical application of the \ZS product will be to combine a
monoid and a group under the hypotheses of Lemma
\ref{PartMulPropsII}\tref{MonoidItem}.  If \(U\) is a monoid and
\(A\) is a group, then the previous paragraph applies and we can
regard \(U\) and \(A\) as substructures of the product.  If Lemma
\ref{PartMulPropsII}\tref{MonoidItem} applies, then there is a
global identity on \(U\times A\).  If we have \(\alpha\in A\) and
\(u\in U\), then we get \(\alpha^u\in A\) which is then invertible.
A formula for the inverse is obtained by looking at
\((\alpha^{-1}\alpha)^u = (\alpha^{-1})^{(\alpha\cdot u)}\alpha^u\).
However, this product must be 1 since
\((\alpha^{-1}\alpha)^u=1^u=1\) by Corollary \ref{HowOneActs}.  Thus
we have \mymargin{TheInvForm} \begin{equation} \label{TheInvForm}
(\alpha^u)^{-1} =(\alpha^{-1})^{(\alpha\cdot u)}.\end{equation}

\subsection{Least common left multiples} We continue our list of
properties on the product with a study of least common left
multiples.  This is the most complicated.

\begin{lemma}\mylabel{PartMulPropsLCLM} Let \((U, D_U, m_U)\) and
\((A, D_A, m_A)\) be sets with partial multiplications and with
mutual actions \((\alpha ,u)\mapsto \alpha \cdot u\) and \((\alpha
,u)\mapsto \alpha ^u\) defined on \(H\subseteq A\times U\).  Assume
\tref{ClosureRules}(a--c) and put the partial multiplication
\tref{TheMult} on \(E = U\times A\).  Assume further the hypotheses
of Lemma \ref{PartMulPropsII}\tref{MonoidItem}.  Then the following
holds.

{\begin{enumerate}\useitem
\renewcommand{\theenumi}{{\itshape \roman{enumi}}}

\item If \(U\) is a cancellative monoid with least common left
multiples, if \(A\) is a group, and if \tref{CoConfRule} holds, then
the multiplication on \(U\times A\) has least common left multiples.
Further, the least common left multiple
\((r,\alpha)(u,\theta)=(s,\beta)(v,\phi)\) of \((u,\theta)\) and
\((v,\phi)\) in \(U\times A\) can be constructed so that
\(r(\alpha\cdot u)=s(\beta\cdot v)\) is the least common left
multiple of \((\alpha\cdot u)\) and \((\beta\cdot v)\) in \(U\).
Lastly, if \(U\times A\) is cancellative (e.g., \tref{CancRules}(a)
holds), then any least common left multiple
\((r,\alpha)(u,\theta)=(s,\beta)(v,\phi)\) of \((u,\theta)\) and
\((v,\phi)\) in \(U\times A\) has the property that \(r(\alpha\cdot
u)=s(\beta\cdot v)\) is the least common left multiple of
\((\alpha\cdot u)\) and \((\beta\cdot v)\) in \(U\).

\end{enumerate}}  \end{lemma}

We need an intermediate lemma.

\begin{lemma}\mylabel{LCLMInt} Under the assumptions of Lemma
\ref{PartMulPropsLCLM}, let \(u\) and \(v\) be in \(U\) with a least
common left multiple \(l=wu=zv\), let \(\gamma\) be in \(A\) and let
\(\alpha=\gamma^w\) and \(\beta=\gamma^z\).  Then \[q=\gamma\cdot l
= (\gamma\cdot q)(\alpha\cdot u) = (\gamma\cdot z)(\beta\cdot v)\]
is the least common left multiple of the pair \((\alpha\cdot u,
\beta\cdot v)\).  \end{lemma}

\begin{proof} Since \(q=\gamma{}\cdot l = \gamma{}\cdot (wu) =
\gamma{}\cdot (zv)\), we have \[\begin{split} q &= (\gamma{}\cdot
w)((\gamma{})^w\cdot u) = (\gamma{}\cdot w)(\alpha{}\cdot u),
\quad\mathrm{and} \\ q &= (\gamma{}\cdot z)((\gamma{})^z\cdot v) =
(\gamma{}\cdot z)(\beta{}\cdot v).\end{split}\]

Thus \(q\) is a common left multiple of \((\alpha{}\cdot u,
\beta{}\cdot v)\) and this pair has a least common left multiple by
hypothesis.  We will use Lemma \ref{MaxesUnique}(c) to argue that
\(q\) is a least common left multiple by showing that any common
left factor of \((\alpha{}\cdot u, \beta{}\cdot v)\) is a unit.

%XXXXXXXXXXXX above and below is a use of common left factor.

Suppose \((\gamma\cdot w)=fh\) and \((\gamma\cdot z)=fj\).  Now
\(w=(\gamma^{-1}\cdot f)((\gamma^{-1})^f\cdot h)\) and
\(z=(\gamma^{-1}\cdot f)((\gamma^{-1})^f\cdot j)\).  This makes
\((\gamma^{-1}\cdot f)\) a common left factor of \((w,z)\) and a
unit by Lemma \ref{MaxesUnique}(c).  Thus \((\gamma^{-1}\cdot
f)k=1\) for some \(k\) and \[1=\gamma\cdot 1 =
\gamma\cdot((\gamma^{-1}\cdot f)k) = f(\gamma^{(\gamma^{-1}\cdot
f)}\cdot k)\] by Corollary \ref{HowOneActs} and \(f\) is a
unit.\end{proof}

\begin{proof}[Proof of Lemma \ref{PartMulPropsLCLM}] We will use the
remarks in Paragraph \ref{MonGpPars} and write \(u\alpha\) instead
of the pair \((u,\alpha)\) in \(U\times A\).

For \(u\) and \(v\) in \(U\) and \(\theta\) and \(\phi\) in \(A\),
we are to find a least common left multiple of \(u\theta\) and
\(v\phi\) assuming that they have a common left multiple.  It is
trivial that this pair has a common left multiple if and only if
\(u\) and \(v\phi\theta^{-1}\) have a common left multiple.
Further, \(p(u\theta)=q(v\phi)\) will be the least common left
multiple of \(u\theta\) and \(v\phi\) if and only if
\(pu=q(v\phi\theta^{-1})\) is the least common left multiple of
\(u\) and \(v\phi\theta^{-1}\).  Thus we shall assume \(u\) and
\(v\) in \(U\) with \(\phi\) in \(A\) and we will look for a least
common left multiple of \(u\) and \(v\phi\) assuming that they have
a common left multiple.

We start with \((x\alpha)(u) = (y\beta)(v\phi)\) or \(x(\alpha\cdot
u)(\alpha^{u}) = y(\beta\cdot v)(\beta^{v}\phi)\). This means that
two equations hold: \(x(\alpha\cdot u) = y(\beta\cdot v)\), and
\(\alpha^{u} = \beta^{v}\phi\).

By assumption there are \(r\) and \(s\) so that \mymargin{ProdLCLMI}
\begin{equation} \label{ProdLCLMI} r(\alpha\cdot u) = s(\beta\cdot
v)\end{equation} is the least common left multiple of the pair
\((\alpha\cdot u, \beta\cdot v)\).

From \(r(\alpha\cdot u) = s(\beta\cdot v)\) \(\alpha^{u} =
\beta^{v}\phi\) we know \(r(\alpha\cdot u)\alpha^{u} = s(\beta\cdot
v)\beta^{v}\phi\) or \mymargin{ProdLCLMII}
\begin{equation}\label{ProdLCLMII}(r\alpha) u = (s\beta)
(v\phi).\end{equation} This is our candidate least common left
multiple for the pair \((u, v\phi)\).  If we show that it is a least
common left multiple, then it satisfies the last sentence of the
statement of Lemma \ref{PartMulPropsLCLM}.

Now assume that \((z\gamma)(u) = (w\delta)(v\phi)\) or
\(z(\gamma\cdot u)(\gamma^{u}) = w(\delta\cdot v)(\delta^{v}\phi)\).
This means that two equations hold: \(z(\gamma\cdot u) =
w(\delta\cdot v)\), and \(\gamma^{u} = \delta^{v}\phi\).  From the
two equalities, \(\alpha^{u} = \beta^{v}\phi\) and \(\gamma^{u} =
\delta^{v}\phi\) we get two expressions for \(\phi\) which must be
equal.  Thus we have \( (\beta^{v})^{-1}(\alpha^{u}) = \phi =
(\delta^{v})^{-1}(\gamma^{u})\) or \( (\delta^{v})(\beta^{v})^{-1} =
(\gamma^{u})(\alpha^{u})^{-1}\).

Using \tref{TheInvForm}, we can work on each side.  The left side
becomes \[(\delta^{v})(\beta^{v})^{-1} =
(\delta^{v})(\beta^{-1})^{(\beta\cdot {v})} =
(\delta)^{(\beta^{-1}\cdot \beta\cdot v)} (\beta^{-1})^{(\beta\cdot
{v})} = (\delta\beta^{-1})^{(\beta\cdot v)}.\] Similarly, the right
side becomes \((\gamma^{u})(\alpha^{u})^{-1} =
(\gamma\alpha^{-1})^{(\alpha\cdot u)}\).

We thus have the equation \[(\gamma\alpha^{-1})^{(\alpha\cdot u)} =
(\delta\beta^{-1})^{(\beta\cdot v)}\] and we have the least common
left multiple \tref{ProdLCLMI} of the pair \((\alpha\cdot u,
\beta\cdot v)\).

The assumption \tref{CoConfRule} and Lemma \ref{StrongCoConf} says
there is an \(\eta\) so that \(\gamma\alpha^{-1} = (\eta)^r\) and
\(\delta\beta^{-1} = (\eta)^s\).  Lemma \ref{LCLMInt} says that \[q
= (\eta\cdot r)((\gamma\alpha^{-1})\cdot (\alpha\cdot u)) =
(\eta\cdot s)( (\delta\beta^{-1})\cdot (\beta\cdot v)) \] is the
least common left multiple of the pair \(((\gamma\alpha^{-1})\cdot
(\alpha\cdot u), (\delta\beta^{-1})\cdot (\beta\cdot v))\).
Equivalently, \[q = (\eta\cdot r)(\gamma\cdot u) = (\eta\cdot
s)(\delta\cdot v)\] is the least common left multiple of the pair
\((\gamma\cdot u, \delta\cdot v)\).

From the equality \(z(\gamma\cdot u) = w(\delta\cdot v)\) we know
that there is a \(k\) so that the common left multiple
\(z(\gamma\cdot u) = w(\delta\cdot v)\) of the pair \((\gamma\cdot
u, \delta\cdot v)\) satisfies \[z(\gamma\cdot u) = w(\delta\cdot v)
= k (\eta\cdot r)(\gamma\cdot u) = k (\eta\cdot s)(\delta\cdot
v). \] Thus \( z = k (\eta\cdot r)\) and \(w = k (\eta\cdot s)\).

We have the candidate least common left multiple \tref{ProdLCLMII}
for the pair \((u, v\phi)\).  We consider \[(k\eta)(r\alpha) u =
(k\eta)(s\beta) (v\phi).\] The left side becomes \[ (k\eta)(r\alpha)
u = k(\eta\cdot r)(\eta)^r \alpha u = z (\gamma\alpha^{-1})\alpha u
= (z\gamma)u\] and the right side becomes \[ (k\eta)(s\beta) (v\phi)
= k(\eta \cdot s)(\eta)^s \beta(v\phi) =
w(\delta\beta^{-1})\beta(v\phi) = (w\delta)(v\phi)\] which is what
we wanted.  

We have found one least common left multiple
\((r,\alpha)(u,\theta)=(s,\beta)(v,\phi)\) of \((u,\theta)\) and
\((v,\phi)\) in \(U\times A\) so that \(r(\alpha\cdot
u)=s(\beta\cdot v)\) is the least common left multiple of
\((\alpha\cdot u)\) and \((\beta\cdot v)\) in \(U\).  If
\tref{CancRules}(a) holds, then \(U\times A\) is cancellative by
Lemma \ref{PartMulProps}\tref{CancItem} and the corresponding
conclusion for left cancellative since we know that the fact that
\(A\) is a group implies that \tref{CancRules}(b) also holds.  Let
\((r',\alpha')(u,\theta)=(s',\beta')(v,\phi)\) be another least
common left multiple of \((u,\theta)\) and \((v,\phi)\).  From Lemma
\ref{MaxesUnique}(a), we know that there is an invertible element
\((t,\gamma)\) so that
\begin{equation}\label{LCLMEqI}(r',\alpha')(u,\theta)=(s',\beta')(v,\phi)
=
(t,\gamma)(r,\alpha)(u,\theta)=(t,\gamma)(s,\beta)(v,\phi).\end{equation}

Expanding and setting the first coordinates equal in the third and
fourth expressions gives
\begin{equation}\label{LCLMEqII}t(\gamma\cdot
r)((\gamma^r\alpha)\cdot u) = t(\gamma\cdot s)((\gamma^s\beta)\cdot
v).\end{equation} From Lemma \ref{LCLMInt}, we know that
\[(\gamma\cdot r)((\gamma^r\alpha)\cdot u) = (\gamma\cdot
s)((\gamma^s\beta)\cdot v)\] is a least common left multiple of
\(((\gamma^r\alpha)\cdot u)\) and \(((\gamma^s\beta)\cdot v)\) in
\(U\).  An easy calculation based on the fact that \((t,\gamma)\) is
invertible shows that \(t\) is invertible in \(U\), so Lemma
\ref{MaxesUnique}(a) shows that \tref{LCLMEqII} is also a least
common left multiple of \(((\gamma^r\alpha)\cdot u)\) and
\(((\gamma^s\beta)\cdot v)\) in \(U\).  From \tref{LCLMEqI} and
cancellativity, we have \((r',\alpha') = (t,\gamma)(r,\alpha)\)
which implies that \(r'=t(\gamma\cdot r)\) and
\(\alpha'=\gamma^r\alpha\).  Similarly, \(s'=t(\gamma\cdot s)\) and
\(\beta'= \gamma^s\beta\).  Combining this with \tref{LCLMEqII}
gives \(r'(\alpha'\cdot u)=s'(\beta'\cdot v)\) is a least common
left mutliple in \(U\) of \((\alpha'\cdot u)\) and \((\beta'\cdot
v)\) which is what we wanted.  \end{proof}

\subsection{Presentations} Because of the difficulties of embedding
the factors in the \ZS product unless there are identities, we
will not attempt to describe presentations of \ZS products of
semigroups.  In all probability, the \ZS product of two semigroups
presented as \(\langle X_i\mid R_i\rangle \), \(i=1,2\), would need
a generating set based on \(X_1\times X_2\).  In the monoid case, we
can get by with \(X_1\cup X_2\).

\begin{lemma}\mylabel{TwistedPres} Under the hypotheses of Lemma
\ref{PartMulPropsII}\tref{MonoidItem}, if the monoids \(U\) and
\(A\) have a presentations \(\langle X \mid R\rangle\) and \(\langle
Y \mid T\rangle\), respectively, with \(X\cap Y=\emptyset\), then a
presentation for the structure defined by \tref{TheMult} on
\(U\times A\) is
\mymargin{ZappaPresI}\begin{equation}\label{ZappaPresI}\langle X\cup
Y \mid R\cup T\cup W\rangle\end{equation} in which \(W\) consists of
all pairs \((\alpha u,(\LeftAct{\alpha}{u})(\alpha^u))\) for
\((\alpha,u)\in A\times U\).  \end{lemma}

\begin{proof} We send \(x\in X\) to \((x,1)\) and \(y\in Y\) to
\((1,y)\).  The comments in Paragraph \ref{MonGpPars} show that this
extends to a homomorphism from the presentation to the \ZS product
which takes any word in \(X\) representing some \(u\in U\) to
\((u,1)\) and any word in \(Y\) representing some \(\alpha\in A\) to
\((1, \alpha)\).  The relations allow us to reduce any word in the
generators to a word of the form \(u\alpha\) where \(u\) is a word
in the alphabet \(X\) and \(\alpha\) is a word in the alphabet
\(A\).  This word is sent to \((u,\alpha)\).  If such a word
\(u\alpha\) is sent to \((1,1)\), then \(u\) is the identity for
\(U\) and must reduce to the empty word using \(R\) and \(\alpha\)
is the identity in \(A\) which must reduce to the empty word using
\(T\).  \end{proof}

Even if the presentations \(\langle X\mid R\rangle\) and \(\langle
Y\mid T\rangle\) of the previous lemma are small, the presentation
\tref{ZappaPresI} is large since \(W\) has an entry for every pair
in \(A\times U\).  Under some rather special hypotheses, we get a
smaller presentation that is sometimes complete.

\begin{lemma}\mylabel{TwistedPresII} Under the hypotheses of Lemma
\ref{PartMulPropsII}\tref{MonoidItem}, if the monoids \(U\) and
\(A\) have a presentations \(\langle X \mid R\rangle\) and \(\langle
Y \mid T\rangle\), respectivley, with \(X\cap Y=\emptyset\), if the
function \((\alpha,u)\mapsto \alpha\cdot u\) takes \(Y\times X\)
into \(X\), and if the function \((\alpha,u)\mapsto \alpha^u\) takes
\(Y\times X\) into \(Y\), then a presentation for the structure
defined by \tref{TheMult} on \(U\times A\) is
\mymargin{ZappaPresII}\begin{equation}\label{ZappaPresII}\langle
X\cup Y \mid R\cup T\cup W\rangle\end{equation} in which \(W\)
consists of all pairs \((\alpha u,(\LeftAct{\alpha}{u})(\alpha^u))\)
for \((\alpha,u)\in Y\times X\).  \end{lemma}

\begin{proof} The proof is word for word identical to the proof of
Lemma \ref{TwistedPres}.  The justifications for the words are
slightly different.  \end{proof}

Question: if the presentations \(\langle X \mid R\rangle\) and
\(\langle Y \mid T\rangle\) in Lemma \ref{TwistedPresII} are
complete, then must the presentation \tref{ZappaPresII} be complete?

The hypotheses of Lemma \ref{TwistedPresII} might be hard to attain.
We discuss some weaker hypotheses.  Assume that presentations
\(\langle X \mid R\rangle\) and \(\langle Y \mid T\rangle\) of the
monoids \(U\) and \(A\), respectively, are given with \(X\cap
Y=\emptyset\), and that functions \(Y\times X\rightarrow Y^*\)
written \((\alpha,u)\mapsto \alpha^u\) and \(Y\times X\rightarrow
X\) written \((\alpha,u)\mapsto \alpha\cdot u\) are given.  The
unequal treatment of the codomains (\(X\) in one case and \(Y^*\) in
the other) is deliberate.

We extend these to functions \(Y^*\times X^*\rightarrow Y^*\) and
\(Y^*\times X^*\rightarrow X^*\) as follows.  Form the monoid
presentation
\mymargin{ZappaPresInt}\begin{equation}\label{ZappaPresInt}\langle
X\cup Y \mid W\rangle\end{equation} in which \(W\) is regarded as a
set of rewriting rules and consists of all pairs \((\alpha
u\rightarrow (\LeftAct{\alpha}{u})(\alpha^u))\) for \((\alpha,u)\in
Y\times X\).

\begin{lemma}\mylabel{ZapPresIntCompl} The presentation
\tref{ZappaPresInt} is complete.  \end{lemma}

\begin{proof} Let \(w\) be a word in the alphabet \(X\cup Y\) and
assume that there are \(n\) appearances of elements of \(X\) in
\(w\) and that they are numbered from left to right in \(w\).  Let
\(C_w=(c_1,c_2,\ldots,c_n)\) where \(c_i\) is the number of elements
of \(Y\) in \(w\) that are to the left of the \(i\)-th element of
\(X\) in \(w\).  An application of \(\rightarrow\) involves a
particular consecutive pair \(\alpha u\) in \(w\) with \(\alpha\in
Y\) and \(u\in X\).  If the pair involves the \(i\)-th element of
\(X\) in \(w\), then applying \(\rightarrow\) keeps all \(c_j\) with
\(j<i\) the same, lowers \(c_i\) by one, and is unpredictable on the
\(c_j\) for \(j>i\).  It is standard that no infinite sequence of
such changes can be applied to an \(n\)-tuple of elements of \(\N\),
so \(\rightarrow\) is terminating.  It is easy to check that
\(\rightarrow\) is strongly coherent since two applications of
\(\rightarrow\) to the same word are either the same or involve
disjoint pairs in \(w\).  \end{proof}

The irreducibles are of the form \(u\alpha\) with \(u\) a word in
the alphabet \(X\) and \(\alpha\) a word in the alphabet \(Y\).
This expresses the monoid presented by \tref{ZappaPresInt} as a \ZS
product of the free monoids \(X^*\) and \(Y^*\).  From Lemma
\ref{MutualActs} we get our desired extensions and in addition get
the fact that they satisfy the four identities
\tref{AssocRules}(a--d).

The extensions to the functions clearly satisfy \tref{FullRule} and
since all structures under discussion are monoids, Corollary
\ref{HowOneActs} implies that they also satisfy
\tref{IdentRules}(a--f).  We are now in a position to state the
following lemma.  Its importance is that the hypotheses are stated
entirely in terms of \(X\), \(Y\), \(R\) and \(T\) without reference
to \(X^*\), \(Y^*\), \(R'\) or \(T'\).  In another paper, the author
uses this lemma to prove that a certain \ZS product used in
combining the Thompson groups with the braid groups has a given
presentation.

\begin{lemma}\mylabel{TwistedPresIII} Assume that presentations
\(\langle X \mid R\rangle\) and \(\langle Y \mid T\rangle\) of the
monoids \(U\) and \(A\), respectively, are given with \(X\cap
Y=\emptyset\), and that functions \(Y\times X\rightarrow Y^*\)
written \((\alpha,u)\mapsto \alpha^u\) and \(Y\times X\rightarrow
X\) written \((\alpha,u)\mapsto \alpha\cdot u\) are given.  Let
\(\sim_R\) and \(\sim_T\) denote the equivalence relations on
\(X^*\) and \(Y^*\), respectively, imposed by the relation sets
\(R\) and \(T\), respectively.

Let the functions be extended to \(Y^*\times X^*\) as above and
assume that they satisfy the following.  If \((u,v)\) is in \(R\),
then for all \(\alpha\in Y\) we have \((\alpha\cdot u, \alpha\cdot
v)\) or \((\alpha\cdot v, \alpha\cdot u)\) is in \(R\) and
\(\alpha^u\sim_T\alpha^v\).  If \((\alpha,\beta)\) is in \(T\), then
for all \(u\in X\) we have \(\alpha\cdot u=\beta\cdot u\) and
\(\alpha^u\sim_T\beta^u\).  Then the extensions induce well defined
functions \(A\times U\rightarrow A\) and \(A\times U\rightarrow U\)
that satisfy the hypotheses of Lemma \ref{TwistedPresII} and thus
the conclusions of that lemma hold.  Further, the restriction of the
function \(A\times U\rightarrow U\) to \(A\times X\) has its image
in \(X\).  \end{lemma}

The assymetry in the hypotheses is deliberate.  Note that the last
sentence implies that the equality \(\alpha\cdot u=\beta\cdot u\)
demanded when \((\alpha,\beta)\) is in \(T\) and \(u\) is in \(X\)
is an equality of elements of \(X\).

\begin{proof} If \(u\) is in \(X\) and \(\alpha\) is in \(Y^*\),
then the identity \((\alpha\beta)\cdot u= \alpha\cdot
(\beta\cdot u)\) and our assumption that \((\alpha,u)\mapsto
\alpha\cdot u\) has image in \(X\) (as opposed to \(X^*\)) implies
that \(\alpha\cdot u\) is in \(X\).  When it is shown that the
hypothesized functions induce well defined functions on \(A\times
U\), we will automatically have the last sentence of the statement.

It suffices to prove that \(\alpha\sim_T\alpha'\)
implies both \((\alpha\cdot u)\sim_R(\alpha'\cdot u)\) and
\(\alpha^u\sim_T(\alpha')^{u}\), and also that \(u\sim_R u'\)
implies both \((\alpha\cdot u)\sim_R(\alpha\cdot u')\) and
\(\alpha^u\sim_T\alpha^{u'}\).  Note that we have restricted
versions of these conclusions as hypotheses of the lemma.  Our task
to extend what is known.

Step 1. Given \((u,u')\) in \(R\) and \(\alpha\in Y^*\), we will
show that \((\alpha\cdot u, \alpha\cdot u')\) or \((\alpha\cdot u',
\alpha\cdot u)\) is in \(R\).  This follows immediately from
\((\alpha\beta)\cdot u= \alpha\cdot (\beta\cdot u)\) and the
assumption that what we want is true when \(\alpha\) is in \(Y\).

Step 2. Given \(\alpha\) and \(\alpha'\) in \(Y^*\) with
\(\alpha\sim_T \alpha'\) and \(u\in X\), we will show that
\(\alpha\cdot u=\alpha'\cdot u\).  It suffices to assume that
\(\alpha\rightarrow_T \alpha'\).  We have \(\alpha=\theta
\beta\rho\) and \(\alpha' = \theta\beta'\rho\) with \((\beta,
\beta')\) in \(T\).  Now \((\theta\beta\rho)\cdot u =
\theta\cdot(\beta\cdot (\rho\cdot u)))\) and
\((\theta\beta'\rho)\cdot u = \theta\cdot(\beta'\cdot (\rho\cdot
u)))\).  Our assumptions are that \(\rho\cdot u\) is some \(v\in
X\), and that \(\beta\cdot v=\beta'\cdot v\).  This completes the
step.

Step 3. Given \(\alpha\) and \(\alpha'\) in \(Y^*\) with
\(\alpha\sim_T \alpha'\) and \(u\in X\), we will show that
\(\alpha^u\sim_T (\alpha')^u\).  It suffices to assume that
\(\alpha\rightarrow_T \alpha'\).  We have \(\alpha=\theta
\beta\rho\) and \(\alpha' = \theta\beta'\rho\) with \((\beta,
\beta')\) in \(T\).  Now \[(\theta\beta\rho)^u =
(\theta^{(\beta\rho)\cdot u})(\beta^{\rho\cdot u})(\rho^u)\] and
\[(\theta\beta'\rho)^u = (\theta^{(\beta'\rho)\cdot
u})((\beta')^{\rho\cdot u})(\rho^u).\] The third factors are
identical.  The common exponent \(\rho\cdot u\) in the second
factors is an element of \(X\) by hypothesis and so the second
factors are related by \(\sim_T\) by hypothesis.  The exponents in
the first factors are \((\beta\rho)\cdot u = \beta\cdot(\rho\cdot
u)\) and \((\beta'\rho)\cdot u = \beta'\cdot(\rho\cdot u)\).  With
\(\rho\cdot u\) in \(X\), our hypotheses say that these exponents
are equal and so the first factors are equal.  The fact that
\(\sim_T\) is a congruence relation completes the step.

Step 4. Given \(\alpha\) and \(\alpha'\) in \(Y^*\) with
\(\alpha\sim_T \alpha'\) and \(u\in X^*\), we will show that
\(\alpha^u\sim_T (\alpha')^u\).  This follows from Step 3 and the
fact that \(\alpha^{(uv)} = (\alpha^u)^v\).

Step 5. Given \(\alpha\) and \(\alpha'\) in \(Y^*\) with
\(\alpha\sim_T \alpha'\) and \(u\in X^*\), we will show that
\(\alpha\cdot u=\alpha'\cdot u\).  If \(u\) is in \(X\), we are done
by Step 2.  If \(u=vw\) with \(v\in X\) and \(w\in X^*\), then
\[\alpha\cdot u=\alpha\cdot(vw)=(\alpha\cdot v)(\alpha^v\cdot w),\]
and \[\alpha'\cdot u=\alpha'\cdot(vw)=(\alpha'\cdot
v)((\alpha')^v\cdot w).\] By Step 2, the first factors are equal.
By Step 3, \(\alpha^v\sim_T(\alpha')^v\) and we are done by
induction on the length of \(u\).

Step 6. Given \((u,u')\) in \(R\) and \(\alpha\in Y^*\), we will
show that \(\alpha^u\sim_T \alpha^{u'}\).  We are done by hypothesis
if \(\alpha\) is in \(Y\).  If \(\alpha=\beta\gamma\) with
\(\beta\in Y^*\) and \(\gamma\in Y\), then \[\alpha^u =
(\beta\gamma)^u = (\beta^{\gamma\cdot u})(\gamma^u),\] and
\[\alpha^{u'} = (\beta\gamma)^{u'} = (\beta^{\gamma\cdot
u'})(\gamma^{u'}).\] We have \(\gamma^u\sim_T\gamma^{u'}\) by
hypothesis.  We also have \((\gamma\cdot u, \gamma\cdot u')\) or
\((\gamma\cdot u', \gamma\cdot u)\) in \(R\) by hypothesis and we
are done by induction on the length of \(\alpha\).

Step 7. Given \(u\) and \(u'\) in \(X^*\) with \(u\sim _R u'\) and
\(\alpha\in Y^*\), we will show that \(\alpha^u\sim_T \alpha^{u'}\).
It suffices to assume \(u\rightarrow_T u'\).  We have \(u=xyz\) and
\(u'=xy'z\) with \((y,y')\) in \(R\).  Now \(\alpha^u =
((\alpha^x)^y)^z\) and \(\alpha^{u'} = ((\alpha^x)^{y'})^z\).
Letting \(\beta=\alpha^x\) in \(Y^*\), we have \(\gamma=\beta^y\)
and \(\gamma'=\beta^{y'}\) related by \(\sim_T\) by Step 6.  Now
\(\gamma^z\sim_T(\gamma')^z\) by Step 4.

Step 8.  Given \(u\) and \(u'\) in \(X^*\) with \(u\sim_R u'\) and
\(\alpha\in Y^*\), we will show that \(\alpha\cdot u\sim_R
\alpha\cdot {u'}\).  It suffices to assume \(u\rightarrow_R u'\).
We have \(u=xyz\) and \(u'=xy'z\) with \((y,y')\) in \(R\).  Now
\[\alpha \cdot u = \alpha\cdot (xyz) = (\alpha\cdot x)(\alpha^x\cdot
y)(\alpha^{(xy)}\cdot z),\] and \[\alpha \cdot u' = \alpha\cdot
(xy'z) = (\alpha\cdot x)(\alpha^x\cdot y')(\alpha^{(xy')}\cdot z).\]
The first factors are identical.  One of the pairs \((\alpha^x\cdot
y, \alpha^x\cdot y')\) or \((\alpha^x\cdot y', \alpha^x\cdot y)\) is
in \(R\) by Step 1.  We have \(\alpha^{(xy)}\sim_T\alpha^{(xy')}\)
by Step 7 and we have that the third factors are equal by Step 5.

Steps 4, 5, 7 and 8 finish showing that the extensions to
\(Y^*\times X^*\) induce well defined functions on \(A\times U\).

Multiplication on \(A\) and \(U\) and the induced functions on
\(A\times U\) are defined by use of representatives.  Since the
multiplications and functions are full, it does not matter which
representatives are chosen.  Thus equations in the monoid are
preserved and \tref{AssocRules}(a--d) must hold for the induced
functions as do \tref{IdentRules}(a--d).  Similarly
\tref{IdentRules}(e,f) hold since the form that they take in monoids
is that the identities are preserved.  Lastly, we get
\tref{FullRule} from the fact that the extensions are defined
on all of \(Y^*\times X^*\).  \end{proof}

\subsection{Associativity of the \ZS product} We would like to show
that the \ZS product is associative in that the map
\((u,(v,a))\mapsto ((u,v),a)\) gives an isomorphism from \(U\zapprod
(V\zapprod A)\) to \((U\zapprod V)\zapprod A\).  It turns out that
it is often easy to establish the isomorphism once it is known that
the two expressions \(U\zapprod (V\zapprod A)\) and \((U\zapprod
V)\zapprod A\) make sense.  What is complicated is making sense out
of the two expressions.  In \ZS products of more than three factors,
there are more expressions to interpret.

Our first step is to assume that we are working with the internal
\ZS product.  This allows us to get isomorphisms from the fact that
two or more expressions accurately represent one fixed
multiplicative structure.

There are still problems with the internal \ZS product.  In the
internal product, the expression \(U\zapprod V\) implies that \(U\)
and \(V\) are substructures of \(U\zapprod V\).  However, the
expression \((U\zapprod V)\zapprod A\) invites the possibility that
the substructure \((U\zapprod V)\) of \((U\zapprod V)\zapprod A\)
contains neither \(U\) nor \(V\).  If one insists that these
containments hold, then some properties must be assumed to make them
so.  We have been unable to find better hypotheses than the direct
assumption that these containments hold.

Note that there does not seem to be a problem if the expressions
\(U\zapprod (V\zapprod A)\) and \((U\zapprod V)\zapprod A\) are both
assumed valid.  After all, if all products are internal, then \(U\)
and \(V\) are in \(U\zapprod V\) and \(V\) and \(A\) are in
\(V\zapprod A\) by definitions.  However, if only one expression,
say \(U\zapprod (V\zapprod A)\), is assumed to be valid, then we
have no reason to assume that \(U\) and \(V\) are contained in
\(UV\) nor can we assume that \(UV\) is even closed under the
product structure.  This gets worse in longer expressions such as
\(U\zapprod (V\zapprod (A\zapprod B))\).

One way out of this is to assume that the containing structure is a
monoid.  Another method, given below, is to make a large number of
containment assumptions.

We start with one approach to the situation involving three factors.
We will modify our approach somewhat when more factors are involved.

We assume the following.  There is a set with partial
multiplication \((M,D,m)\) and three subsets \(U\), \(V\) and \(A\)
of \(M\) that satisfy the six statements below.

\begin{enumerate}

\item All of \(U\), \(V\) and \(A\) are closed under \(m\) in that
if \(xDy\) and \(x\) and \(y\) are in \(U\) (for example), then
\(xy\in U\).

\item The two products \(UV\) and \(VA\) are ``unique factorization
products'' in that if \(y\in UV\) (for example), then \(y\) is
uniquely expressible as \(y=uv\) with \(u\in U\) and \(v\in V\).

\item The two subsets \(UV\) and \(VA\) are ``semi-permutable'' in
that \(VU\sub UV\) and \(AV\sub VA\).

\item The containments \(U\sub UV\), \(V\sub UV\cap VA\) and \(A\sub
VA\) all hold.

\item The set with mutliplication \((M,D,m)\) is categorical.

\item Every element of \(M\) is uniquely expressible as \(uw\) with
\(u\in U\) and \(w\in VA\).

\end{enumerate}

The last assumption is the only one that does not treat \(U\) and
\(A\) symmetrically.

The assumptions (1) and (3) and the associativity implicit in (5)
easily give that \(UV\) and \(VA\) are closed under the
multiplication.  The categorical property is inherited by both
\(UV\) and \(VA\).  From (1), (2), (4), (5) and (6) we can write
\(M\simeq U\zapprod (V\zapprod A)\) in which the isomorphism takes
\(uva\) to \((u,(v,a))\), and we know that every \(x\in M\) is
uniquely expressible as \(x=u(va)\) with \(u\in U\), \(v\in V\) and
\(a\in A\).  Now every \(x\in M\) is expressible as \(ya\) with
\(y\in UV\) and \(a\in A\), and \(uva\mapsto ((u,v),a)\) is an
isomorphism from \(M\) to \((U\zapprod V)\zapprod A\).  Thus we have
proven the following.

\begin{lemma}\mylabel{ZappaAssoc} Given a set with partial
multiplication \((M,D,m)\) and subsets \(U\), \(V\) and \(A\) of
\(M\) satisfying (1--6) above, the multiplication gives mutual
actions for all of the pairs \((U,V)\), \((V,A)\), \((UV,A)\) and
\((U,VA)\) and the maps \((u,(v,a))\mapsto uva \mapsto ((u,v),a)\)
give isomorphisms \(U\zapprod (V\zapprod A)\rightarrow M\rightarrow
(U\zapprod V)\zapprod A\). \end{lemma}

As remarked above, Lemma \ref{ZappaAssoc} holds if \(M\) is a
monoid.

Abused notation can mislead.  Given \(a\in A\) and \(u\in U\), we
have \(au=u'v'a'\) with \(u'\in U\), \(v'\in V\) and \(a'\in A\).
It is tempting to write \(au=(a\cdot u)(a^u)\) and make conclusions.
In fact from the structure of \(M\) as \((U\zapprod V)\zapprod A\),
we get \(au=(a\cdot u)_1(a^u)_1\) with \((a\cdot u)_1\in UV\) and
\((a^u)_1\in A\) so that \((a\cdot u)_1=u'v'\) and \((a^u)_1=a'\).
From the structure of \(M\) as \(U\zapprod (V\zapprod A)\), we get
\(au=(a\cdot u)_2(a^u)_2\) with \((a\cdot u)_2\in U\) and
\((a^u)_2\in VA\) so that \((a\cdot u)_2=u'\) and \((a^u)_2=v'a'\).

To discuss higher iterations of the \ZS product, we need some
notation.

If \(\odot\) is a binary operation we let
\mymargin{ParenOne}\begin{equation}\label{ParenOne}E(X_1\odot
X_2\odot X_3\odot \cdots \odot X_{n-1}\odot X_n)\end{equation}
denote a parenthesization so that the order of operations is well
defined.  If the \(X_i\) are replaced by \(Y_i\) and \(\odot\) is
replaced by another binary operation \(\boxdot\), then
\mymargin{ParenTwo}\begin{equation}\label{ParenTwo}E(Y_1\boxdot
Y_2\boxdot Y_3\boxdot \cdots \boxdot Y_{n-1}\boxdot
Y_n)\end{equation} using the same symbol \(E\) will represent the
``same'' parenthesization of the expression as in \tref{ParenOne}.
That is, if the parenthesization scheme imposed by \(E\) is
represented by a binary tree, then the same binary tree is used to
control the parenthesization in both \tref{ParenOne} and
\tref{ParenTwo}.  We will also use the lowly comma ``\(,\)'' as a
binary operation that creates ordered pairs.  With \(E\) the same
symbol as in \tref{ParenOne} and \tref{ParenTwo}, we get
\[E(x_1,x_2,x_3, \ldots, x_{n-1}, x_n)\] as a well defined
parenthesization that reduces the ordered \(n\)-tuple \((x_1,
\ldots, x_n)\) to an heirarchy of ordered pairs based on the same
binary tree used in \tref{ParenOne} and \tref{ParenTwo}.

We now consider a set with partial multiplication \((M,D,m)\) and a
finite number of subsets \(U_1, U_2, \ldots, U_n\) of \(M\).  For
\(1\le i\le j\le n\) we let \(W_i^j=U_iU_{i+1}\cdots U_j\).  Note
that \(W_i^i=U_i\).  We make the following assumptions.

{\begin{enumerate} \renewcommand{\theenumi}{\alph{enumi}}

\item The set with partial multiplication \((M,D,m)\) is
categorical.

\item Each \(x\in M\) is uniquely expressible as \(x=u_1u_2\cdots
u_n\) with each \(u_i\in U_i\).

\item For \(i<n\) and \(u_i\in U_i\) there is a \(u_{i+1}\in
U_{i+1}\) so that \(u_iu_{i+1}\) is defined and equal to \(u_i\),
and for \(1<j\) and \(u_j\in U_j\) there is a \(u_{j-1}\in U_{j-1}\)
so that \(u_{j-1}u_j\) is defined and equal to \(u_j\).

\item For \(1\le i\le j\le n\), the set \(W_i^j\) is closed under
the multiplication.

\end{enumerate}}

It follows from (a) that all the \(W_i^j\) are categorical.  For
\(i\le k\le l\le j\), it follows from (a) and (c) that \(W_k^l\sub
W_i^j\).  Thus for \(1\le i \le j \le n\), we have \(W_i^j\sub M\)
and we get from (b) that every \(w\in W_i^j\) is uniquely
expressible as \(w=u_iu_{i+1}\cdots u_j\) with each \(u_k\in U_k\).
For \(1\le i\le j< k\le n\), Lemma \ref{ZappaRecon} says that
sending \((u_iu_{i+1}\cdots u_j, u_{j+1}\cdots u_k)\) to
\(u_iu_{i+1}\cdots u_k\) gives an isomorphism from \(W_i^j\zapprod
W_{j+1}^k\) to \(W_i^k\).  We now inductively have the following.

\begin{lemma} Let \((M,D,m)\) be a set with partial multiplication,
let \(U_1, U_2, \ldots, U_n\) be a finite number of subsets of
\(M\), and assume that (a--d) above hold.  Then for each
parenthesization \(E(U_1\zapprod U_2\zapprod \cdots \zapprod U_n)\),
sending \(E(u_1, u_2, u_3, \ldots ,u_{n-1}, u_n)\) to \(u_1u_2
\cdots u_n\) is an isomorphism from \(E(U_1\zapprod U_2\zapprod
\cdots \zapprod U_n)\) to \(M\). \end{lemma}

%\bibliography{thompson}
\providecommand{\bysame}{\leavevmode\hbox to3em{\hrulefill}\thinspace}

\noindent Department of Mathematical Sciences

\noindent State University of New York at Binghamton

\noindent Binghamton, NY 13902-6000

\noindent USA

\noindent email: matt@math.binghamton.edu

\end{document}